\RequirePackage{fix-cm}
\documentclass[smallextended]{svjour3}       
\smartqed  
\usepackage{CJK}
\usepackage{graphicx}
\usepackage{epstopdf}
\usepackage{amsmath}
\usepackage{amsfonts,amssymb}
\usepackage{rotating}
\usepackage{subfigure}
\usepackage{multirow}
\usepackage{dcolumn}
\usepackage{booktabs}
\usepackage{color}
\usepackage{caption}
\usepackage{float}
\usepackage{longtable}
\usepackage{cite}
\usepackage{array}
\renewcommand{\Re}{{\rm I}\!  {\rm R}}
\newcommand{\be}{\begin{equation}}
\newcommand{\ee}{\end{equation}}
\newcommand{\bit}{\begin{itemize}}
\newcommand{\eit}{\end{itemize}}
\newcommand{\sgn}{\hbox{sgn}}
\newtheorem{algorithm}{Algorithm}
\begin{document}

\title{A Semismooth Newton Method for Support Vector Classification and Regression
}


\author{Juan Yin         \and
        Qingna Li 
}


\institute{Juan Yin \at
              School of Mathematics and Statistics, Beijing Institute of Technology, Beijing, 100081, China\\
              \email{2120171435@bit.edu.cn}     
           \and
           Qingna Li \at
              Corresponding author. School of Mathematics and Statistics, Beijing Key Laboratory on MCAACI, Beijing Institute of Technology,
		Beijing, 100081, China\\
   \email{qnl@bit.edu.cn}. This author's research was supported by the National Science Foundation of China(No.11671036).
}

\date{Received: date / Accepted: date}

\maketitle

\begin{abstract}
Support vector machine is an important and fundamental technique in machine learning. In this paper, we apply a semismooth Newton method to solve two typical SVM models: the L2-loss SVC model and the $\epsilon$-L2-loss SVR model. The semismooth Newton method is widely used in optimization community. A common belief on the semismooth Newton method is its fast convergence rate as well as high computational complexity. Our contribution in this paper is that by exploring the sparse structure of the models, we significantly reduce the computational complexity, meanwhile keeping the quadratic convergence rate. Extensive numerical experiments demonstrate the outstanding performance of the semismooth Newton method, especially for problems with huge size of sample data (for \texttt{news20.binary} problem with 19996 features and 1355191 samples, it only takes three seconds). In  particular, for the $\epsilon$-L2-loss SVR model, the semismooth Newton method significantly outperforms the leading solvers including DCD and TRON. 
\end{abstract}

\section{Introduction}
\label{intro}
Support vector machine (SVM) is a popular and important statistical learning technique {\cite{Q1, Q9, Q10, Q11, Q12, Lee}}. SVMs hold records in performance benchmarks for handwritten digit recognition, text categorization, information retrieval, and time-series prediction. They  are commonly used in the analysis of DNA micro-array data \cite{Q6, Q7, Q8, Q9, Q10}. Two main categories for support vector machines (SVMs) are support vector classification (SVC) and support vector regression (SVR). Support vector classification is a learning machine for two-group classification problems \cite{Q16}. The support vector regression was extended from SVC by Boser et al. \cite{Q23}.   Most of the optimization methods for SVM models solve the dual problems, partly due to some nonsmooth properties of the primal functions. Two typical examples are the L2-loss SVC model and the $\epsilon$-L2-loss SVR model. Below we give a brief review on methods for the above two models, which motivate the work in our paper. For a survey of optimization methods for machine learning, we refer to \cite{largesvm2017,Friedman2001}.

 For the L2-loss SVC model, due to the nondifferentiability of the gradient of the objective function, Mangasarian \cite{Q5} introduced a finite Newton method. It is basically a semismooth Newton method with unit step size, and the inverse of Hessian matrix is used to calculate the Newton direction. Keerthi and DeCoste \cite{Q17} proposed a modified Newton method. They compute the Newton point and do an exact line search to determine step length. A trust region Newton method (TRON) \cite{Q4} was proposed for the L2-loss SVC model. Chang et al. \cite{Q13} proposed a coordinate descent method for the  primal problem and Hsieh et al. \cite{Q3} proposed a dual coordinate descent method (DCD) for the dual problem  of the L2-loss SVC model.   Very recently, Hsia et al. \cite{Q18}\footnote{We realized this work when we drafted our paper.} performed a study on trust region update rules in Newton's  method. 
For the $\epsilon$-L2-loss SVR model,  Ho and Lin \cite{Q14} applied the TRON and DCD to solve it and a smoothing Newton method was proposed by Gu et al. \cite{Q15} .  To deal with large scale of data, stochastic gradient methods become popular in solving large scale SVM models \cite{largesvm2017}. Stochastic gradient method and its variants have good performance in machine learning \cite{largesvm2017}. Classical stochastic gradient descent (SGD) was proposed by Robbins and Monro in 1951 \cite{SGD}. Johnson and Zhang \cite{SVRG} proposed an accelerating stochastic gradient descent using predictive variance reduction (SVRG). Recently, Tan et al. \cite{SVRG-BB} put forward to use the Barzilai-Borwein (BB) method to automatically compute step sizes for SGD and SVRG, which leads to two algorithms: SGD-BB and SVRG-BB. In their paper, numerical results suggest that SVRG-BB and SGD-BB clearly outperform SVRG and SGD respectively. To summarize, one can see that despite the competitiveness of Newton-type methods in SVM,  little attention has been paid to the semismooth Newton method in solving the two models.

 On the other hand, in optimization community, the semismooth Newton method has been well studied, and has been successfully used in many applications, especially in solving modern optimization problems, such as the nearest
correlation matrix problem \cite{QiSun2006,Qi2013A}, the nearest Euclidean distance matrix
problem \cite{LiQi2012}, the tensor eigenvalue complementarity problem \cite{chen}, solving the system of absolute value equations \cite{Cruz}, the solution of quasi-variational inequations \cite{Facchinei}, as well as linear and convex quadratic semidefinite programming problems \cite{Zhao2009}. The concept of semismoothness was introduced by Mifflin \cite{Mifflin}, and was popularized by Qi and Sun \cite{QiSun93}. In \cite{QiSun93}, Qi and Sun proposed  a nonsmooth version of the classical Newton's method. Compared with the classical Newton method, the semismooth Newton method can solve nonsmooth equations, meanwhile can keep the local quadratic convergence rate under certain conditions. A semismooth Newton method  was extended to solve the nonsmooth matrix equations by Qi and Sun \cite{QiSun2006}. Recently, the semismooth Newton method has been frequently used to solve some important problems, for example, Lasso problems \cite{LiSun2017}, OSCAR and SLOPE models \cite{Sun2017}, approximating
weighted time series of finite rank \cite{Qi2017}, and convex clustering problems\cite{YuanSun2018}.

 Compared with the wide usage of the semismooth Newton method in optimization community, little attention has been paid to the semismooth Newton method in machine learning, especially in SVM models. In this paper, we will set up such a bridge by applying a globalized semismooth Newton method to models of SVC and SVR, i.e., the L2-loss SVC model and the $\epsilon$-L2-loss SVR model. 
 A common belief on the semismooth Newton method is its fast convergence rate as well as high computational complexity. Our contribution in this paper is that by exploring the sparse structure of the models, we significantly bring down the computational complexity, meanwhile keeping the quadratic convergence rate. Another advantage is that it is able to handle  the case with a huge number of sample data, since it solves the primal problem rather than the dual. Extensive numerical experiments demonstrate the outstanding performance of the semismooth Newton method, especially for problems with huge size of sample data (for \texttt{news20.binary} problem with 19996 features and 1355191 samples, it only takes about three seconds). In  particular, for the $\epsilon$-L2-loss SVR model, the semismooth Newton method significantly outperforms the leading solvers including DCD and TRON.


The remaining parts of this paper are organized as follows. In Section \ref{sec-model}, we introduce the formulation of two models of SVMs, i.e., the L2-loss SVC model and the $\epsilon$-L2-loss SVR model. In Section \ref{sec-algorithm}, we introduce the semismooth Newton method and apply it to solve the two mentioned models. In Section \ref{sec-sparsity}, we characterize the generalized Jacobian of the objective functions in the two models, and highlight how to maintain the quadratic convergence rate and reduce the computational complexity by making use of the sparse structure. In Section \ref{sec-numerical}, we collect test data from LIBLINEAR, a popular package for SVM,  and conduct extensive numerical experiments to show the efficiency of our algorithm.  We also  do   comparisons with other state-of-art solvers, such as TRON, DCD and SVRG-BB. Finally, we conclude our paper in Section \ref{sec-conclusions}.

\section{Two Models of SVMs}
\label{sec-model}
{\bf The L2-loss SVC Model}
Given training data $x_1,x_2,\dots,x_l\in\Re^n$ and the corresponding label $ y_1,y_2,\dots,y_l\in\{-1,1\}$, the basic idea of support vector classification is to find a hyperplane $\omega^Tx+b=0$ to separate the data, where $\omega\in\Re^n$ and $b\in\Re$ are unknown parameters. The traditional SVM model is
\be\label{svm}
\begin{split}
\min_{\omega\in\Re^n,\ b\in \Re} & \frac12\|\omega\|_2^2  \\
\hbox{s.t.}\ \ \ \ \ &  y_i(\omega^Tx_i+b)\ge 1,\ i= 1,\dots,\ l.
\end{split}
\ee
Here we actually require that the data should be strictly separated, i.e., the constraints must be satisfied strictly.
This model is based on the assumption that the data can be linearly separated. In practice, one usually using the following regularized model which
 allows that the data could be wrongly labelled, i.e., the inequality constraints can be violated
\be\label{hing}
\min_{\omega\in\Re^n,\ b\in\Re} \ \frac12\|\omega\|^2 + C\sum_{i = 1}^l \xi(\omega; x_i, y_i, b)
\ee
where $C>0$ is a penalty parameter and $\xi(\cdot)$ is the loss function. If $\xi(\omega; x_i, y_i, b)=\max(1-y_i(\omega^Tx_i+b), 0)$, it is referred as the $l_1$ hinge loss function; if $\xi(\omega; x_i, y_i, b)=\max(1-y_i(\omega^Tx_i+b), 0)^2 $, it is the squared hinge-loss function  which we call  L2-loss function; if $\xi(\omega; x_i, y_i, b)=\log(1+e^{-y_i(\omega^Tx_i+b)})$, it is referred as logistic loss function. In our paper, we focus on the L2-loss SVC model, i.e.,
\be\label{prob-svc}
\min_{\omega\in\Re^n, b\in\Re} \ \frac12\|\omega\|^2 + C\sum_{i = 1}^l \max(1-y_i(\omega^Tx_i+b), 0)^2.
\ee
Recent works on support vector classification often omit the bias term because it hardly affects the performance on most data \cite{Q3,Q14}. Therefore, by appending each instance with an additional dimension:
\be
x_i^T\leftarrow[x_i^T, 1], \ \ \ \ \ \omega^T\leftarrow[\omega^T, b],
\ee
we obtain the following model, which is  the first model we will consider (referred as L2-loss SVC \cite{Q3}).
\be\label{prob}
\min_{\omega\in\Re^n} f_1(\omega):=\frac12\|\omega\|^2+C\sum_{i=1}^l \max(1-y_i(\omega^Tx_i), 0)^2.\\
\ee

{\bf The $\epsilon$-L2-loss SVR Model}
Given training data $x_1,x_2,\dots,x_l\in\Re^n$ and the corresponding  observations $ y_1,y_2,\dots,y_l$,
SVR is to find $\omega\in\Re^{n}$ such that $\omega^Tx_i+b$ is close to the target value $y_i$, $i = 1,\dots, l$. The $\epsilon$-L2-loss SVR model (Similarly, we omit the bias term $b$ for SVR) is as follows
\be\label{svr}
\min_{\omega\in\Re^n}f_2(\omega):=\frac12\|\omega\|^2+C\sum_{i=1}^l \xi_\epsilon(\omega; x_i, y_i), \\
\ee
where
\be\label{loss2}
\xi_\epsilon(\omega; x_i, y_i)=\max(|\omega^Tx_i-y_i|-\epsilon, 0)^2\\
\ee
is the $\epsilon$-insensitive loss function which we call $\epsilon$-L2-loss function associated with $(x_i, y_i)$. The parameter $\epsilon>0$ is given so that the loss is zero if $|\omega^Tx_i-y_i|\leq \epsilon$. Ho and Lin \cite{Q14}, and Gu et al. \cite{Q15} refer to SVR using $(\ref{loss2})$ as L2-loss SVR and $\epsilon$-SVR respectively. We refer to it as  $\epsilon$-L2-loss SVR.

\begin{figure}[H]
  \centering
  \begin{minipage}[t]{.45\linewidth}
\includegraphics[width=1\textwidth]{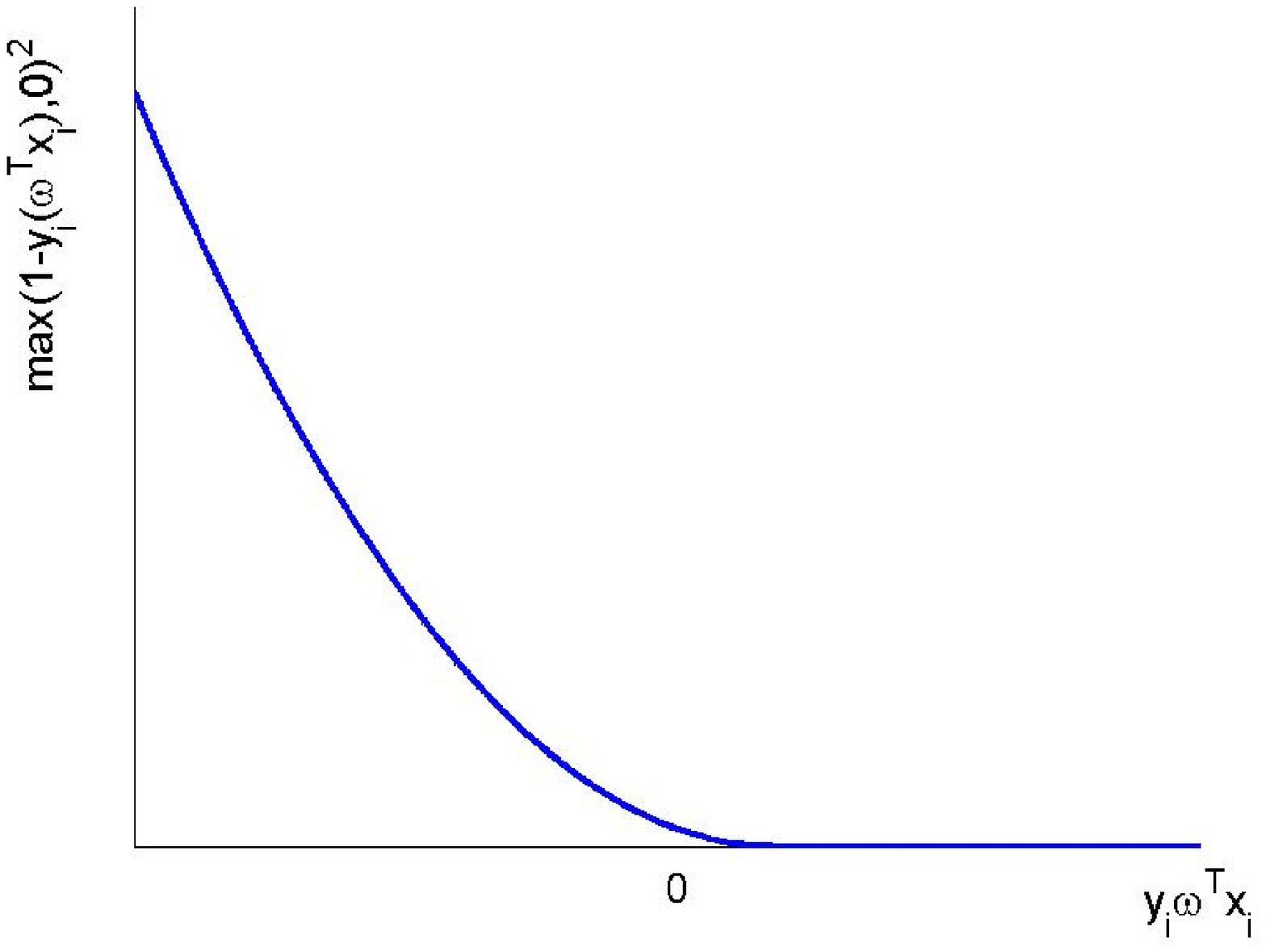}
  \caption*{(a)\ \ L2-loss function for SVC}
  \end{minipage}
  \begin{minipage}[t]{.45\linewidth}
\includegraphics[width=1\textwidth]{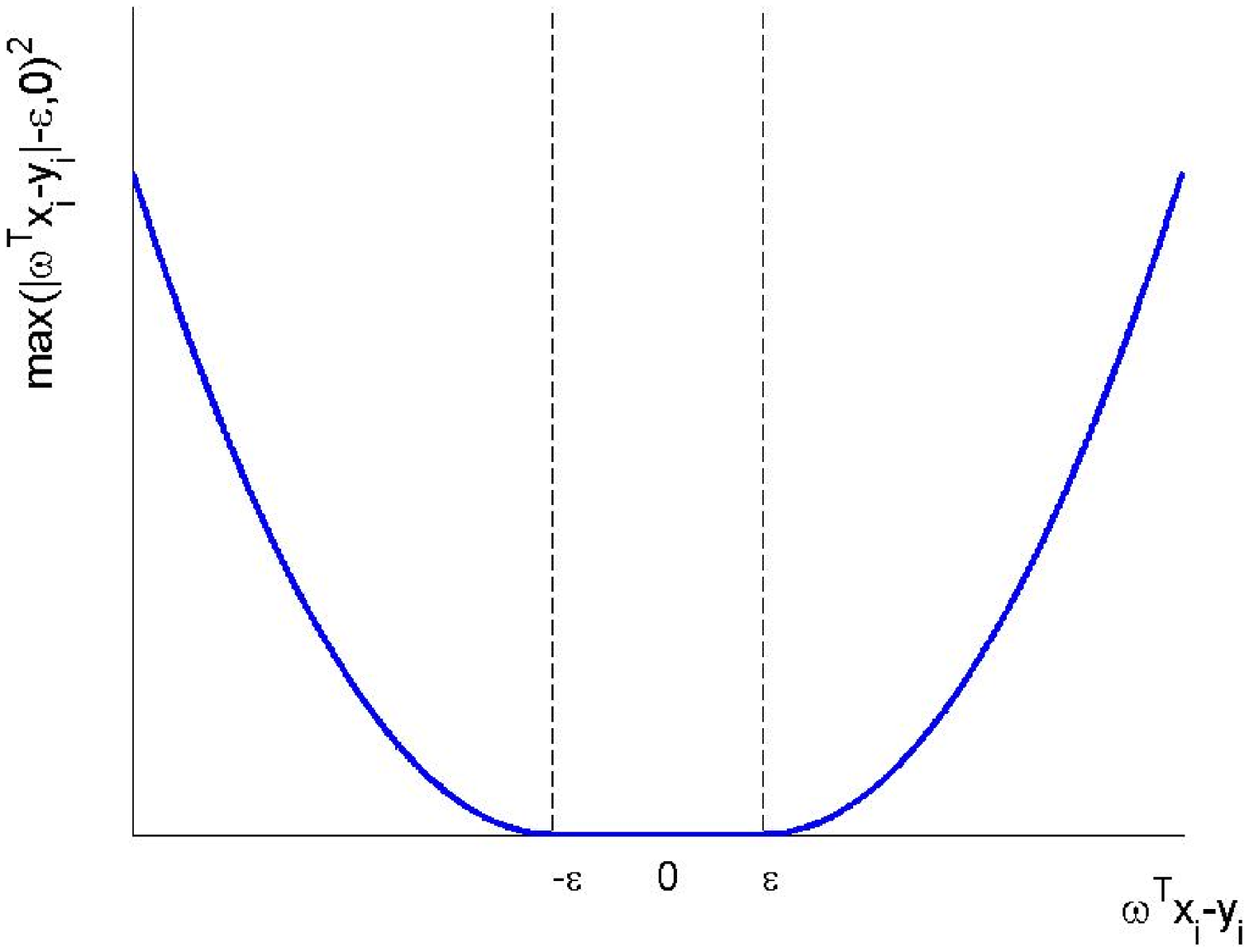}
  \caption*{(b)\ \ $\epsilon$-L2-loss function for SVR}
  \end{minipage}
  \caption{Demonstration of Loss Functions}\label{L2-loss functions}
\end{figure}

One can easily verify that   the functions of  (\ref{prob})  and  (\ref{loss2}) are continuously differentiable but not twice differentiable. An illustration of the loss functions is in Figure $\ref{L2-loss functions}$.
\section{A Semismooth Newton Method}\label{sec-algorithm}

In this section, we will apply the semismooth Newton method to solve (\ref{prob})  and  (\ref{svr}). It is divided into two parts. In the first part, we introduce some preliminaries. In the second part, we apply the semismooth Newton method to solve (\ref{prob})  and  (\ref{svr}). 
\subsection{Preliminaries}
\label{sec:2}
In this part, we will introduce some preliminaries about the semismooth Newton method. The semismoothness of a function is closely related to the generalized Jacobian in the sense of Clarke \cite{Clarke}, which is stated as follows.

Let $\Phi: \Re^m\to \Re^l$ be a (locally) Lipschitz function. According to Rademacher's theorem \cite[Sect. 14]{Redemacher}, $\Phi$ is differentiable almost everywhere. Define
\[
D_{\Phi}:=\{x\in\Re^m\ |\ \Phi \hbox{ is differentiable at } x\}.
\]
Let $\Phi'(x)$ denote the Jacobian of $\Phi$ at $x\in D_\Phi$. The Bouligand subdifferential of $\Phi$ at $x\in\Re^m$ is then defined by
\[
\partial_B\Phi(x):=\{V\in\Re^{m\times l}\ |\ V\hbox{ is an accumulation point of }\Phi'(x^k), \ x^k\to x, \ x^k\in D_\Phi\}.
\]
The generalized Jacobian in the sense of Clarke \cite{Clarke}  is the convex hull of $\partial_B\Phi(x)$, i.e.,
 \[
\partial\Phi(x) = \hbox{co} \ \partial_B \Phi(x),
\]
where $\hbox{co} (\partial_B \Phi(x))$ is the convex hull of $\partial_B \Phi(x)$.
The concept of semismoothness was introduced by Mifflin \cite{Mifflin} for functionals. It was extended to vector-valued functions by Qi and Sun \cite{QiSun93}.

\begin{definition}
We say that $\Phi$ is semismooth at $x$ if (i) $\Phi$ is directional differentiable at $x$ and (ii) for any $V\in\partial \Phi(x+h)$,
\[
\Phi(x+h) -\Phi(x) -Vh = o(\|h\|).
\]
$\Phi$ is said to be strongly semismooth at $x$ if $\Phi$ is semismooth at $x$ and for any $V\in\partial \Phi(x+h)$,
\[
\Phi(x+h) -\Phi(x) -Vh = O(\|h\|^2).
\]
\end{definition}

Some particular examples for semismooth functions are as follows. \bit
\item Piecewise linear functions are strongly semismooth.
\item The composition of (strongly) semismooth functions is also (strongly) semismooth.
\eit

For example, according to the definition above, $\max (0, t)$ is strongly semismooth and  the Clarkes' generalized gradient of $\max (0, t)$ is
\be\label{partial}
\partial \max(0,\ t) =\left \{
\begin{array}{ll}
1,\  & t>0,\\
0,\ & t<0,\\
 v,\ 0\le v\le1, & t=0.
\end{array}
\right.
\ee
\subsection{A Semismooth Newton Method Applied to (\ref{prob})  and  (\ref{svr}) }

For $\Phi:\ \Re^m\rightarrow \Re^m$, a nonsmooth version of the classical Newton method to solve the equations $\Phi(x)=0$ is as follows
\be\label{seminewton}
x^{k+1} = x^k - V_k^{-1}\Phi(x^k), \ V_k\in\partial \Phi(x^k), \ k = 0,\ 1,\ 2,\  \ldots,
\ee
where $x^0$ is an initial point.  In general, the above iterative method does not converge. However, Qi and Sun \cite{QiSun93} show that if $\Phi$ is semismooth, then the iterate sequence converges superlinearly. It is from then that the semismooth Newton method became popular. We would also like to highlight that if $\Phi$ is continuously differentiable, then $\partial \Phi(x)$ reduces to a singleton, which is the Jacobian of $\Phi(x)$. In this situation, the algorithm is the classical Newton method.


For solving the following problem
\be\label{minf}
\min_{\omega\in\Re^n} f(\omega)
\ee
where $f: = f_1(\omega)$ or $f = f_2(\omega)$. It is easy to verify that $f$ is strongly convex and continuously differentiable with \[\nabla f_1(\omega) = \omega-2C\sum_{i=1}^l \max(0, 1-y_i\omega^Tx_i)y_ix_i \]
and
\[
\nabla f_2(\omega)=\omega+2C\sum_{i=1}^l \max(|\omega^Tx_i-y_i|-\epsilon, 0)\sgn(\omega^Tx_i-y_i)x_i,
\]
where $\sgn(t)$ is defined as $1$ if $t\ge 0$ and $-1$ otherwise. Therefore, solving (\ref{minf}) is equivalent to solving
\be
\nabla f(\omega)=0.
\ee

One can see that $\nabla f_1(\omega)$ and $\nabla f_2(\omega)$ are continuous, but not differentiable. Fortunately, based on our analysis in Section 3.1, we can see that $\nabla f_1(\omega)$ and $\nabla f_2(\omega)$ are strongly semismooth. Therefore, we can apply the semismooth Newton method to solve (\ref{prob}) and (\ref{svr}). In practice, we use the following
 well studied globalized version of the semismooth Newton method \cite[Algorithm 5.1]{QiSun2006}. 

\begin{algorithm}\label{algorithm2}  A globalized semismooth Newton method
\begin{itemize}

\item [S0] Given $k:=0$. Choose $\omega^0,\ \sigma\in(0,\ 1)$, $\rho\in(0,\ 1)$,\ $\delta>0$, and $\eta_0>0,\ \eta_1>0 $.
\item [S1] Calculate $\nabla f(\omega^k)$. If $\|\nabla f(\omega^k)\|\le \delta$, stop. Otherwise, go to S2.
\item [S2] Select an element $V^k\in \partial^2 f(\omega^k)$ and apply CG \cite{Hestenes} to find an approximate solution $d^k$ by
\be\label{linearsystem}
V^k d^k + \nabla f(\omega^k)=0
\ee
such that
\[
\|V^k d^k + \nabla f(\omega^k)\|\le \mu_k \|\nabla f(\omega^k)\|
\]
where  $\mu_k = \min(\eta_0    ,\ \eta_1\|\nabla f(\omega^k)\|)$.

\item [S3] Do line search, and let $m_k>0$ be the smallest integer such that the following holds
\[
f(\omega^k+\rho^md^k)\le f(\omega^k) + \sigma \rho^m\nabla f(\omega^k)^Td^k.
\]
 Let $\alpha_k = \rho^{m_k}$.
 \item [S4] Let $\omega^{k+1} = \omega^k + \alpha_k d^k${\color{red},} $k:=k+1$, go to S1.

\end{itemize}
\end{algorithm}

{\bf Remark.} Note that Mangsarian \cite{Q5} proposed a finite Armijo Newton method for solving L2-loss SVC. Different from Mangsarian's algorithm, we use the conjugate gradient (CG) method proposed by Hestenes and Stiefel \cite{Hestenes} to solve the linear system in S2 for obtaining the descent direction $d$. We note that Hsia et al. \cite{Q18} proposed using line search and trust region to obtain step length but they focused on investigating the trust region update rules in  Newton's method for L2-loss SVC.
\section {Quadratic Convergence Rate and Low Computational Complexity} \label{sec-sparsity}
The tradition view about the Newton method is the fast convergence rate and its expensive computational cost due to the usage of second order information. In this section, we will show that when  the semismooth Newton method is applied to solve the two models (\ref{prob}) and (\ref{svr}), the quadratic convergence rate can be well maintained. Furthermore, we can also reduce the computational complexity dramatically by fully exploring the sparse structure of the models. We divide this section into three parts.
In the first part, we  characterize the generalized Jacobian of $\nabla f(\omega)$, which is used in Alg. \ref{algorithm2}. Then we discuss how to maintain the quadratic convergence rate of the semismooth Newton method. Finally, we will bring down the computational complexity by exploring the sparse structure of the models.

\subsection{Characterization of Generalized Jacobian}
In Algorithm \ref{algorithm2}, we need to calculate $\partial^2 f(\omega)$, i.e., the generalized Jacobian of $\nabla f(\omega)$.
For the L2-loss SVC model (\ref{prob}),   by the chain rule \cite[Theorem 2.3.9]{Clarke},  there is $\partial^2  f_1(\omega)\subseteq  \mathcal V_1$ where
\[
\mathcal V_1 = \{  I + 2C\sum_{i = 1}^l h_ix_ix_i^T,\ h_i\in\partial \max(0,z_i(\omega)), \ z_i(\omega) = 1-y_i\omega^Tx_i, \ i = 1,\  \ldots,\  l \}.
\]

For $\epsilon$-L2-loss SVR, the generalized Jacobian of $\nabla f_2(\omega)$ is characterized  in the following proposition.
\begin{proposition} For $f_2(\omega)$ defined as in (\ref{svr}), there is $\partial^2  f_2(\omega)\subseteq  \mathcal V_2$, where
\[
\mathcal V_2 = \{  I + 2C\sum_{i = 1}^l h_ix_ix_i^T, h_i\in\partial \max(0,z_i(\omega)), \ z_i(\omega) = |\omega^Tx_i-y_i|-\epsilon, \ i = 1, \ldots, l \}.
\]
\end{proposition}
{\bf Proof.}
Recall that
\begin{eqnarray*}
\nabla f_2(\omega)&=&\omega+\left\{
\begin{array}{ll}
2C\sum\limits_{i=1}^l \max(\omega^Tx_i-y_i-\epsilon, 0)x_i,\ & \omega^Tx_i-y_i \geq 0\\
2C\sum\limits_{i=1}^l \max(y_i-\omega^Tx_i-\epsilon, 0)(-x_i),\ & \omega^Tx_i-y_i < 0\\
\end{array}
\right.\\
&=&\omega+2C\sum_{i=1}^l \max(|\omega^Tx_i-y_i|-\epsilon, 0)\sgn(\omega^Tx_i-y_i)x_i.
\end{eqnarray*}
In the following, we first discuss the generalized Jacobian of $\max(|\omega^Tx_i-y_i|-\epsilon, 0)\sgn(\omega^Tx_i-y_i)x_i$, $i = 1, \dots, l$. First, denote $Q_i: \Re^n\rightarrow \Re^n$ by
\[Q_i(\omega)=\max(|\omega^Tx_i-y_i|-\epsilon, 0)\sgn(\omega^Tx_i-y_i)x_i, \ i = 1, \dots, l.\]
There is
\begin{eqnarray*}
Q_i(\omega)
&=&\left\{
\begin{array}{ll}
(\omega^Tx_i-y_i-\epsilon, 0)x_i,\ & \omega^Tx_i-y_i \geq 0, \ \ \hbox{and} \ \  \omega^Tx_i-y_i-\epsilon> 0\\
0,\ & \omega^Tx_i-y_i \geq 0, \ \ \hbox{and} \ \  \omega^Tx_i-y_i-\epsilon= 0\\
0,\ & \omega^Tx_i-y_i \geq 0,\ \ \hbox{and} \ \   \omega^Tx_i-y_i-\epsilon<0\\
-(y_i-\omega^Tx_i-\epsilon, 0)x_i,\ \  & \omega^Tx_i-y_i < 0,\ \ \hbox{and} \ \  y_i-\omega^Tx_i-\epsilon> 0\\
0,\ & \omega^Tx_i-y_i < 0,\ \ \hbox{and}\ \   y_i-\omega^Tx_i-\epsilon=0\\
0,\ & \omega^Tx_i-y_i< 0,\ \ \hbox{and} \ \   y_i-\omega^Tx_i-\epsilon<0\\
\end{array}
\right.\\
&=&\left\{
\begin{array}{ll}
(\omega^Tx_i-y_i-\epsilon, 0)x_i,\ & \omega^Tx_i-y_i > \epsilon\\
0,\ & -\epsilon \leq \omega^Tx_i-y_i \leq\epsilon\\
-(y_i-\omega^Tx_i-\epsilon, 0)x_i,\ \  & \omega^Tx_i-y_i <-\epsilon.\\
\end{array}
\right.\\
\end{eqnarray*}
We can see that
$Q_i(\omega)$ is differentiable when $\omega^Tx_i-y_i >\epsilon$, or $-\epsilon < \omega^Tx_i-y_i <\epsilon$, or $\omega^Tx_i-y_i <-\epsilon$. However, $Q_i(\omega)$ is  not differentiable if $|\omega^Tx_i-y_i|=\epsilon$.
When  $\omega$ satisfies $|\omega^Tx_i-y_i| -\epsilon>0$ or $|\omega^Tx_i-y_i| -\epsilon<0$, the Jacobian of $Q_i(\omega)$  can be easily calculated by
\begin{eqnarray*}
 Q_i'(\omega)&=&\left\{
\begin{array}{ll}
x_ix_i^T,\ & \omega^Tx_i-y_i > \epsilon\\
0,\ & -\epsilon < \omega^Tx_i-y_i<\epsilon\\
x_ix_i^T,\ \  & \omega^Tx_i-y_i <-\epsilon\\
\end{array}
\right.\\
&=&\left\{
\begin{array}{ll}
x_ix_i^T,\ \  & |\omega^Tx_i-y_i| > \epsilon \\
0,\ & |\omega^Tx_i-y_i|<\epsilon\\
\end{array}
\right.\\
\end{eqnarray*}
Next, we calculate the generalized Jacobian of $Q_i(\omega)$ when $\omega^Tx_i-y_i =-\epsilon$ and $\omega^Tx_i-y_i =\epsilon$.
By Section \ref{sec:2}, $D_{Q_i}=\{\omega\in\Re^n\ |\ |\omega^Tx_i-y_i|-\epsilon\neq0\}$, and
\[
\partial_B  Q_i(\omega)=\{V\ |\ V\hbox{ is an accumulation point of } Q_i'(\omega^j), \omega^j\in D_{Q_i}, \omega^j\rightarrow \omega\}.
\]
Consider at $\omega\in\Re^n$ where $\omega^Tx_i-y_i =-\epsilon$. We choose a sequence $\{\omega^j\}\subset D_{Q_i}$, such that $\omega^j\rightarrow \omega$ and $(\omega^j)^Tx_i-y_i=-\epsilon-\frac1j$. Then $\lim\limits_{j\rightarrow\infty}  Q_i'(\omega^j)=x_ix_i^T.$
Similarly, choose another sequence $\{\omega^j\}\subset D_{Q_i}$, such that $(\omega^j)^Tx_i-y_i=-\epsilon+\frac1j$. There is $\lim\limits_{j\rightarrow\infty}  Q_i'(\omega^j)=0$. Then we have: at $\omega\in\Re^n$ with $\omega^Tx_i-y_i =-\epsilon$,
\[
\partial_BQ_i(\omega)=\{0, x_ix_i^T\}.
\]
Consequently, at $\omega\in\Re^n$ with $\omega^Tx_i-y_i =-\epsilon$,
\[
\partial(Q_i(\omega))=\hbox{co}(\partial_BQ_i(\omega))=\hbox{co}(\{0, x_ix_i^T\}).
\]
Similarly, at $\omega\in\Re^n$ with $\omega^Tx_i-y_i =\epsilon$, we have
\[
\partial(Q_i(\omega))=\hbox{co}(\{0, x_ix_i^T\}).
\]
To sum up, we get
\begin{eqnarray*}
\partial Q_i(\omega)&=&\left\{
\begin{array}{ll}
x_ix_i^T,\ \  & |\omega^Tx_i-y_i| > \epsilon, \\
0,\ & |\omega^Tx_i-y_i|<\epsilon,\\
{\hbox{co}(\{0, x_ix_i^T\})}, \ \  & |\omega^Tx_i-y_i| =\epsilon. \\
\end{array}
\right.\\
\end{eqnarray*}
Note that \begin{eqnarray*}
{\hbox{co}(\{0, x_ix_i^T\})}&=&
\{(1-h_i)\cdot0+h_ix_ix_i^T, h_i\in[0,1]\} =
 \{h_ix_ix_i^T, h_i\in[0,1]\}.
\end{eqnarray*}
The generalized Jacobian of $Q_i(\omega)$ is then given by
\begin{eqnarray*}
\partial Q_i(\omega)&=&\left\{
\begin{array}{ll}
x_ix_i^T,\ \  & |\omega^Tx_i-y_i| > \epsilon, \\
0,\ & |\omega^Tx_i-y_i|<\epsilon,\\
h_ix_ix_i^T, h_i\in[0,1], \ \  & |\omega^Tx_i-y_i| =\epsilon \\
\end{array}
\right.\\
&=& \{h_ix_ix_i^T, h_i = 1  \hbox{ if }   |\omega^Tx_i-y_i|> \epsilon ; \ h_i=0 \ \hbox{if}\   |\omega^Tx_i-y_i|<\epsilon;\\
& & \hbox{and} \ h_i\in[0,1] \ \hbox{if}\ |\omega^Tx_i-y_i| = \epsilon\}.
\end{eqnarray*}

 By \cite[Proposition 2.3.3]{Clarke}, we know  that $\partial(\sum Q_i)(\omega)\subset\sum\partial Q_i(\omega)$ if $Q_i(i=1,\ldots, l)$ is a family of functions each of which is Lipschitz near $\omega$. Therefore, we have
\begin{eqnarray*}
\partial^2f_2(\omega)&=&\partial (\nabla f_2(\omega))\\
&\subseteq &I+2C\sum_{i=1}^l\partial Q_i(\omega)\\
&=&\{I+2C\sum\limits_{i=1}^lh_ix_ix_i^T, h_i=1 \ \hbox{if}\  |\omega^Tx_i-y_i|> \epsilon  ; \ \hbox{and} \ h_i=0 \ \hbox{if}\\&&   |\omega^Tx_i-y_i|<\epsilon; \ \hbox{and} \ h_i\in[0,1] \ \hbox{if}\ |\omega^Tx_i-y_i| = \epsilon, \ i = 1, \dots, l \}.
\end{eqnarray*}


By letting $z_i(\omega) = |\omega^Tx_i-y_i|-\epsilon$, and recall (\ref{partial}), we get $\partial^2  f_2(\omega)\subseteq  \mathcal V_2$, where
\[
\mathcal V_2 = \{  I + 2C\sum_{i = 1}^l h_ix_ix_i^T, h_i\in\partial \max(0,z_i(\omega)), \ z_i(\omega) = |\omega^Tx_i-y_i|-\epsilon, \ i = 1, \ldots, l \}.
\] The proof is finished. \hfill$\Box$

\subsection{Local Quadratic Convergence Rate}

The local convergence result for the semismooth Newton method (\ref{seminewton}) is given as follows.

\begin{theorem}\cite[Thm.3.2]{QiSun93} \label{thm-1} let $x^*$
be  a solution of $\Phi(x)=0$ and Let $\Phi$  be  a locally
Lipschitz function which is semismooth  at $x^*$. Assume that   all $V\in \partial \Phi(x^*)$ are nonsingular. Then every sequence generated by  (\ref{seminewton})  is superlinearly convergent to $x^*$, provided that the starting point $x^0$ is sufficiently close to $x^*$. Moreover, if 
 $\Phi$ is strongly semismooth at $x^*$,  the convergence rate is quadratic.

\end{theorem}

From Theorem \ref{thm-1},  to guarantee the local quadratic convergence rate of   Alg. \ref{algorithm2},  we need to check the positive definiteness of each element in $\partial ^2f(\omega)$.
From the characterization of $\partial^2 f_1(\omega)$ and $\partial ^2f_2(\omega)$, one can easily see that for any $V\in \mathcal V_1\bigcup\mathcal V_2$, there is $V\succeq I$. In other words, $V-I$ is positive semidefinite for any $V\in \mathcal V_1 $ and $\mathcal V_2$. Consequently, we have the following proposition.
\begin{proposition}
For any $V\in \mathcal V_i$, $i = 1, 2$,   $V$ is positive definite.
\end{proposition}
Due to the positive definiteness of $V\in\mathcal V_1$ and $V\in\mathcal V_2$,  the local convergence result in Theorem \ref{thm-1} holds, and the semismooth Newton method applied to solve (\ref{prob}) and (\ref{svr}) enjoys local quadratic convergence rate. 

{\bf Remark.} Here we would like to highlight that not only the quadratic convergence of the semismooth Newton method can be guaranteed in theoretical point of view, it can also be verified from the  numerical point of view. In fact, in our numerical test, quadratical convergence rate can always be observed. More details of the quadratic convergence rate are demonstrated in Section 5.1.

\subsection {Exploring Sparsity to Reduce Computational Complexity}

As mentioned before, the traditional view about the semismooth Newton method is its high computational complexity since it needs to calculate the generalized Jacobian. Also, it needs to solve the linear system in order to get the Newton direction. In this part, we will demonstrate our view about the semismooth Newton method. That is, by exploring the sparse structure, the computational cost can be significantly reduced, which is even lower than calculating the Jacobian. Specifically, the computational complexity can be reduced from $O(ln^2)$ to $O(n|I_k|)$, where $|I_k|\ll l$.

We take the L2-loss SVC model as an example. In each iteration $k$, one needs to solve the linear system (\ref{linearsystem}) to get $d^k$. 
 In our algorithm, we solve the linear system by CG, the computational burden then lies in calculation $V\Delta \omega$ for $\Delta \omega\in\Re^n$ in each CG iteration. Below, we will compare the computational cost of calculating $V\Delta\omega$ by traditional implementation and that by our implementation.

 {\bf Traditional Implementation.}

The traditional implementation of calculating $V\Delta\omega$ is to first generate $V^k$ and save it, then calculate $V^k\Delta\omega$. The computational cost in each step is shown in Table \ref{tab-1} where only multiplication and division are taken into account.

\begin{table}
\centering
\caption{Computational Cost for Traditional Implementation}\label{tab-1}
\begin{tabular}{c|c|c}

\hline
 Step & Formula & Computational Cost
\\\hline

Form $V$&$I + 2C\sum_{i = 1}^l h_ix_ix_i^T$ & $l(n^2+n)+n^2+1$\\\hline

Calculate $V\Delta\omega$&$V\Delta\omega$& $n^2$\\\hline
\end{tabular}
\end{table}
The computational complexity for traditional implementation is then $O(ln^2)$.

{\bf Our Implementation.}

  In our implementation, we didn't store $V$ explicitly. Instead, we calculate $V\Delta \omega$ directly by the right hand side of the following formula
\[
{ V\Delta\omega = \Delta\omega+ 2C\sum\limits_{i = 1}^l h_i(x_i^T\Delta\omega) x_i.}
\]
As we can see, one only need to calculate the second term in the right hand {side} of the above formula. Here we would like to highlight that by taking the product of $V\Delta\omega$ directly, we get avoid of forming matrix $x_ix_i^T$. Instead, we can first take the vector product $x_i^T\Delta\omega$, which will result in a scale, then conduct scale-vector multiplication $(x_i^T\Delta\omega)\cdot x_i$. This will lead to the computational cost of $(2n+1)l+n+1$.

Moreover,  recall that $h_i^k\in\partial \max(0,z_i^k)$ and some of the $h_i$'s are actually zero due to the definition in (\ref{partial}). Consequently, for those indices with $h_i=0$, it is not necessary to calculate the item $ h_i(x_i^T\Delta\omega) x_i$. Consequently,  let $I_k:\ =\{i:\ z_i^k(\omega^k) = 1-y_ix_i^T\omega^k >0\}$. At iteration $k$, we choose $h^k\in\Re^l$ in the following way:
\[
h^k_i=\left\{
\begin{array}{ll}
1, & i \in I_k\\
0, & z_i^k <0\\
0,& z_i^k =0,
\end{array}
\right.
\]
As a result, ${ V\Delta\omega}$ will reduce  to
\[{ V\Delta\omega = \Delta\omega + 2C\sum\limits_{i \in I_k}(x_i^T\Delta\omega)  x_i.}\]
The computational cost then becomes $2|I_k|n+n+1$. This is the implementation we use in our code. These are summarized in Table \ref{tab-2}.

\begin{table}
\centering
\caption{Computational Cost for Our Implementation}\label{tab-2}
\begin{tabular}{c|c|c}

\hline
   & Formula & Computational Cost
\\\hline


Calculate $V\Delta\omega$ directly &$V\Delta\omega = \Delta\omega+ 2C\sum\limits_{i = 1}^l h_i(x_i^T\Delta\omega) x_i$&  $l(2n+1)+n+1$ \\\hline
With definition of $I_k$&$ V\Delta\omega = \Delta\omega + 2C\sum\limits_{i \in I_k}(x_i^T\Delta\omega)  x_i$& $2|I_k|n+n+1$\\\hline
\end{tabular}
\end{table}

To further see the size of $I_k$,
note that for $z_i>0$, it means that the $i$-th sample can not be linearly separated, i.e., it violates $1-y_ix_i^T\omega<0$, so we need to penalized the violation. In this case, we actually assume that only few number of such $i$ will happen. Therefore, it means that near the  optimal solution $\omega^*$, $|I_k|$ is much smaller than the sample size $l,$ i.e., $|I_k|\ll l$. We can see that
compared to calculate $V\Delta\omega$ directly, the complexity in each iteration is reduced from $O(ln^2)$ to $O(|I_k|n)$.



In a word, due to the special sparse structure of problem (\ref{prob}) and (\ref{svr}), our way of calculating $V\Delta\omega$ will lead to low computational cost, which is much lower than the classical Newton and semismooth Newton method.

\section{Numerical Results}\label{sec-numerical}

In this section, we analyze the performance of the semismooth Newton method for solving L2-loss SVC and $\epsilon$-L2-loss SVR. It is divided into five parts. In the first part, we demonstrate the low complexity of the semismooth Newton method as well as the quadratic convergence rate. In the second and third parts, we discuss the performance of the semismooth Newton method for L2-loss SVC and $\epsilon$-L2-loss SVR, respectively, due to different choices of parameters. In the fourth part, we compare our algorithm with the methods in LIBLINEAR \cite{Q25}, including trust region Newton method (TRON) and dual coordinate descent method (DCD). In the last part, we compare with SVRG-BB \cite{SVRG-BB}, one of the most efficient stochastic gradient methods.

All experiments are tested in Matlab R2013b in Windows 7 on a Lenovo desktop computer with an Intel(R) Core(TM) i5-3470M CPU at 3.20 GHZ and 4 GB of RAM.
 Throughout the computational experiments, we use the following parameters in the semismooth Newton method:
 $
 \sigma=10^{-4},\  \rho=0.5, \  \omega^0=ones(n,1),\  \delta =10^{-6}, \ {\eta_0=0.05},\  \eta_1=0.5.
$
When solving the linear system by CG, we set the maxium  number iterations as 200.

Due to the different criteria of error evaluation for L2-loss SVC and $\epsilon$-L2-loss SVR, we use the standard real data sets from LIBSVM for classification and regression (42 data sets for classification and 12 data sets for regression). For some datasets of classification whose labels don't belong to $\{-1,1\}$, we change their labels and set them belong to $\{-1,1\}$. For example, for the dataset: $``$breast-cancer$"$, samples' labels are either 2 or 4. We turn the label 2 into -1 and the label 4 into 1. Similarly, we use the same strategy for datasets: $``$liver-disorders$"$, $``$mushrooms$"$, $``$phishing$"$, $``$skin$\_$nonskin$"$ and $``$svmguide1$"$.   
Detailed information of data sets for classification and regression is given in Table $\ref{A1}$ and Table $\ref{A}$.


\begin{table}[H]
\centering
\caption{Data Information for Classification (l is the number of instances, n is the number of features, $\sharp$ nonzeros represents the number of non-zero elements in all training instances and density shows the ratio: $\sharp$ nonzeros/(l$\cdot$n)).}\label{A1}
\begin{tabular}{ccccc}
\hline\noalign{\smallskip}
Data set& l & n & $\sharp$ nonzeros & density\\
\noalign{\smallskip}\hline\noalign{\smallskip}
  a1a & 30956 & 123 & 429343 & 11.28$\%$\\
  a2a & 30296 & 123 & 420188 & 11.28$\%$ \\
  a3a& 29376 & 123 & 407430 & 11.28$\%$ \\
  a4a& 27780 & 123 &385302 & 11.28$\%$\\
  a5a &26147 & 123 & 362653 & 11.28$\%$\\
  a6a& 21341 & 123 &295984 & 11.28$\%$\\
  a7a & 16461 &123 & 228288 & 11.28$\%$\\
  a8a& 22696 &123 & 314815& 11.28$\%$\\
  a9a& 32561 & 123 & 451592 & 11.28$\%$\\
  australian & 690 &14 & 8447 & 87.44$\%$\\
  breast-cancer&638 &10 & 6380 & 100$\%$\\
  cod-rna& 59535 & 8 & 476280& 100$\%$\\
  colon-cancer & 62 & 2000 & 124000 & 100$\%$\\
  diabetes & 768 & 8 & 6135 & 99.85$\%$\\
  duke breast-cancer & 38 & 7129 & 270902 & 100$\%$\\
  fourclass & 862 &2 & 1717& 99.59$\%$\\
  german.numer& 1000 & 24 & 23001 & 95.84$\%$\\
  gisette & 6000 & 5000 & 29729997 & 99.10$\%$\\
  heart & 270 & 13 & 3510 & 100$\%$\\
  ijcnn1 & 49990 & 22 & 649870 & 59.09$\%$\\
  \noalign{\smallskip}\hline
\end{tabular}
\end{table}

\begin{table}[H]
\centering
\begin{tabular}{ccccc}
\hline\noalign{\smallskip}
Data set& l & n & $\sharp$ nonzeros & density\\
\noalign{\smallskip}\hline\noalign{\smallskip}
  ionosphere & 351 & 34& 10551 & 88.41$\%$\\
  leukemia & 38 & 7129 & 270902 & 100$\%$\\
  liver-disorders & 145 & 5 & 725 & 100$\%$\\
  mushrooms & 8124 & 112 & 170604 & 18.75$\%$\\
  news20.binary & 19996 & 1355191 & 9097916 & 0.03$\%$\\
  phishing & 11055 & 68 & 331610 & 44.11$\%$\\
  rcv1.binary & 20242 & 47236 & 1498952 & 0.16$\%$\\
  real-sim & 72309 & 20958 & 3709083 & 0.24$\%$\\
  skin$\_$nonskin &245057 & 3 & 735171 & 100$\%$\\
  splice & 2175 & 60 & 130500 & 100$\%$\\
  sonar & 208 & 60 & 12479 & 99.99$\%$\\
  svmguide1 & 3089 & 4 &12356 & 100$\%$\\
  svmguide3 & 1243 & 22 & 27208 & 99.50$\%$\\
  w1a& 47272& 300 & 551176 & 3.89$\%$ \\
  w2a & 46279 & 300 & 539213 & 3.89$\%$ \\
  w3a & 44837 & 300 & 522338  & 3.89$\%$  \\
  w4a&42383 & 300 & 493583 & 3.89$\%$ \\
  w5a & 39861 & 300 & 464466 & 3.89$\%$ \\
  w6a & 32561 & 300 & 379116 & 3.89$\%$ \\
  w7a& 25057 & 300 & 291438 & 3.89$\%$ \\
  w8a& 49749 & 300 & 579586 & 3.89$\%$ \\
  covtype.binary&581012&54&6940438 & 22.12$\%$ \\
 \noalign{\smallskip}\hline
\end{tabular}
\end{table}
\begin{table}[H]
\centering
\caption{Data Information  for Regression. }\label{A}
\begin{tabular}{cccccc}
\hline\noalign{\smallskip}
Data set& l & n &  $\sharp$ nonzeros & density & range of y\\
\noalign{\smallskip}\hline\noalign{\smallskip}
 abalone & 4177 & 8 & 32080 & 96.00$\%$ &[4, 29] \\
  bodyfat & 252 & 14 & 3528 & 100$\%$ &[1.00, 1.11]\\
   cpusmall & 8192 & 12 & 98304 & 100$\%$ &[0, 99]\\
   tfidf.train&16087&150360&19971015& 0.83$\%$ &[-7.90, -0.52]\\
   tfidf.test & 3308 & 150360 & 4559533 & 0.92$\%$ &[-7.14, -1.69]\\
  eunite2001 & 336 & 16 & 2651& 49.31$\%$ &[612, 876]\\
   housing & 506 & 13 & 6578 & 100$\%$ &[5, 50]  \\
   mg & 1385 & 6 & 8310 & 100$\%$ &[0.42, 1.32] \\
   mpg & 392 & 7 & 2614 & 95.26$\%$ &[9, 46.6] \\
  pyrim & 74 & 27 & 1720 & 86.09$\%$ &[0.1, 0.9] \\
   space$\_$ga & 3107 & 6 & 18642 & 100$\%$ &[-3.06, 0.10] \\
  triazines & 186 & 60 & 9982& 89.44$\%$ &[0.1, 0.9] \\
\noalign{\smallskip}\hline
\end{tabular}
\end{table}
To see the performance of the semismooth Newton method, we  report the following information: the number of iterations $k$, the total number of CG iterations $cg$, the cputime $t$ in second, as well as the final $\|\nabla f(\omega^k)\|$, denoted as $res$. We also use an index of accuracy to further evaluate the quality of the solution returned by  our method.
For L2-loss SVC, let $\hat x_i$ be a test data, 
 the predicted label is then calculated  as follows
\[
\hat y_i=\sgn(\omega^T\hat x_i),
\]
where $\omega$ is generated by the semismooth Newton algorithm.
 The accuracy is then calculated by
 \[
 \frac{\hbox{number of test data whose predicted labels are correct}}{\hbox{number of test set}}.
 \]
For $\epsilon$-L2-loss SVR, we let $y_i'=\omega^T\hat x_i,  i=1,2,\ldots,m $, m is the total number of testing data, $\hat x_i$ is the element of testing data. We use the mean squared error (MSE) to show our algorithm's test accuracy, which is calculated by
\begin{equation}
\hbox{MSE}=\frac1l\sum_{i=1}^m(y_i-y_i')^2,
\end{equation}
where $y_i', $ are the observed data corresponding to $\hat x_i$, $ i =1, \cdots, m$.

\subsection{Demonstration of Low Computational Complexity and Quadratic Convergence Rate}

{\bf Demonstration of Low Computational Complexity and Sparsity.}

As analyzed above, the model of problem (\ref{prob}) we solved has good sparsity. In this part, we will give an example for the description of sparsity. Here we set $C=\frac1l\times 10^2$ for convenience. For dataset: ``covtype.binary", we can see that the semismooth Newton method takes 11 iterations until terminating successfully and the data set contains 581012 instances. Recall $I_k:\ =\{i:\ z_i^k(\omega^k) = 1-y_ix_i^T\omega^k >0\}$. For each iteration, $|I_k|$ is recorded as follows.
\[
|I_k|=[0, 0, 443992,  278855, 148650, 86290, 57933, 47991, 46343, 46292, 46293].
\]
The corresponding ratio of $|I_k|$ over sample size $l$ is calculated by
\[
\hbox{ratio}_k = \frac{|I_k|}{l}=[0, 0, 0.764, 0.480, 0.256, 0.149, 0.100, 0.083, 0.080, 0.080, 0.080]
\]

\begin{figure}[H]
  \centering
  \begin{minipage}[t]{.48\linewidth}
  \includegraphics[width=1\textwidth]{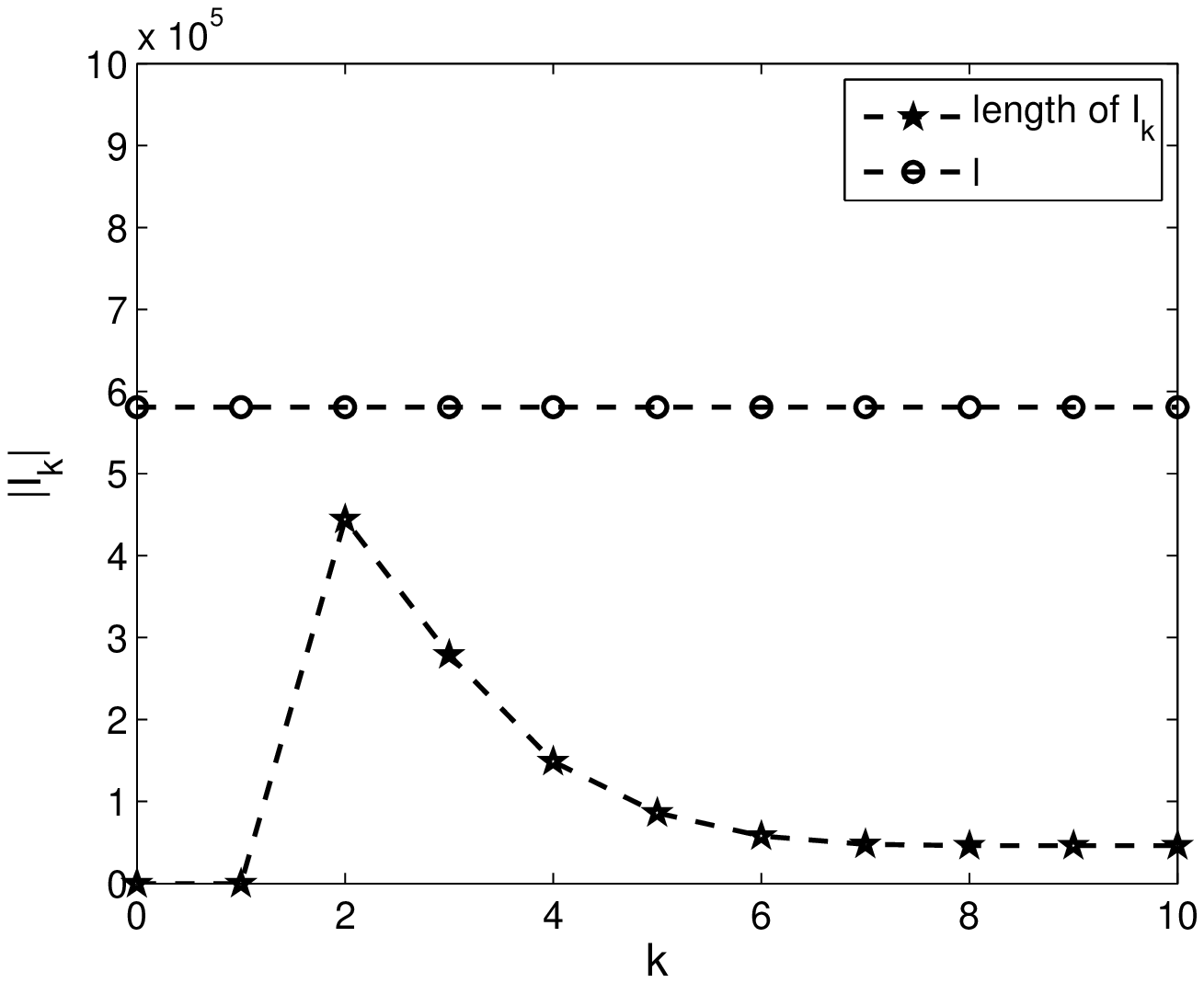}
  \caption*{(a)\ \ The change of $|I_k|$}
  \end{minipage}
  \begin{minipage}[t]{.48\linewidth}
\includegraphics[width=1\textwidth]{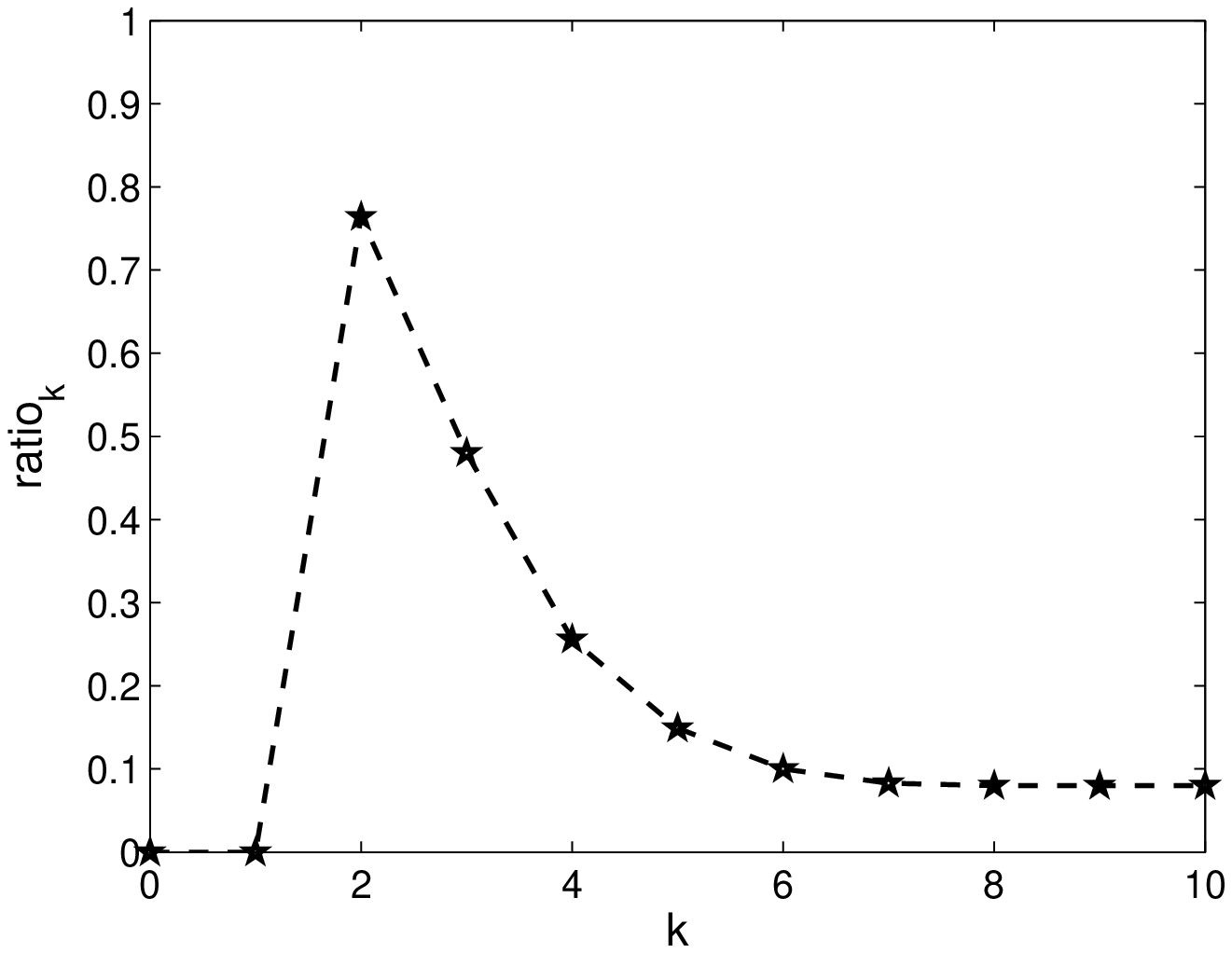}
  \caption*{(b)\ \ The change of $\hbox{ratio}_k$}
  \end{minipage}
  \caption{ Demonstration of $|I_k|$ and $\hbox{ratio}_k$ for Each Iteration}\label{f}
\end{figure}

 We plot $|I_k|$ and $\hbox{ratio}_k$ in Figure \ref{f}. We can see that  $|I_k|$ is always under the horizontal line and $\hbox{ratio}_k$ is always less than 1. In particular, $|I_k|$ is significantly smaller than the total number of instances except $k = 2$ and the value of $\hbox{ratio}_k$ is even less than 0.1 at some iterations, indicating that the computational cost is significantly saved from $O(ln^2)$ to $O(|I_k|n)$.

{\bf Demonstration of Quadratic Convergence Rate.}

\begin{figure}[H]
  \centering
  \begin{minipage}[t]{.48\linewidth}
\includegraphics[width=1\textwidth]{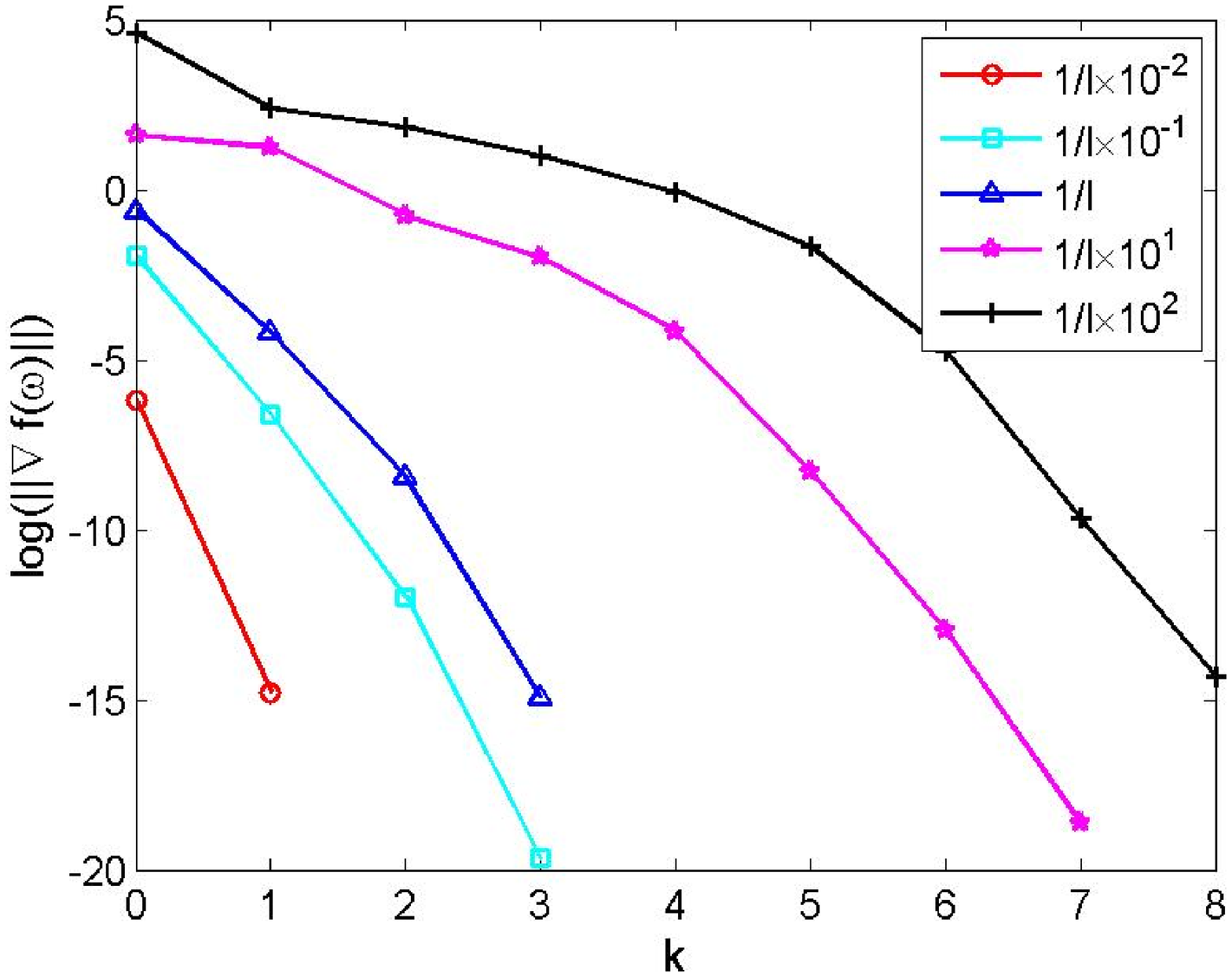}
  \caption*{(a)\ \ w3a}
  \end{minipage}
  \begin{minipage}[t]{.48\linewidth}
\includegraphics[width=1\textwidth]{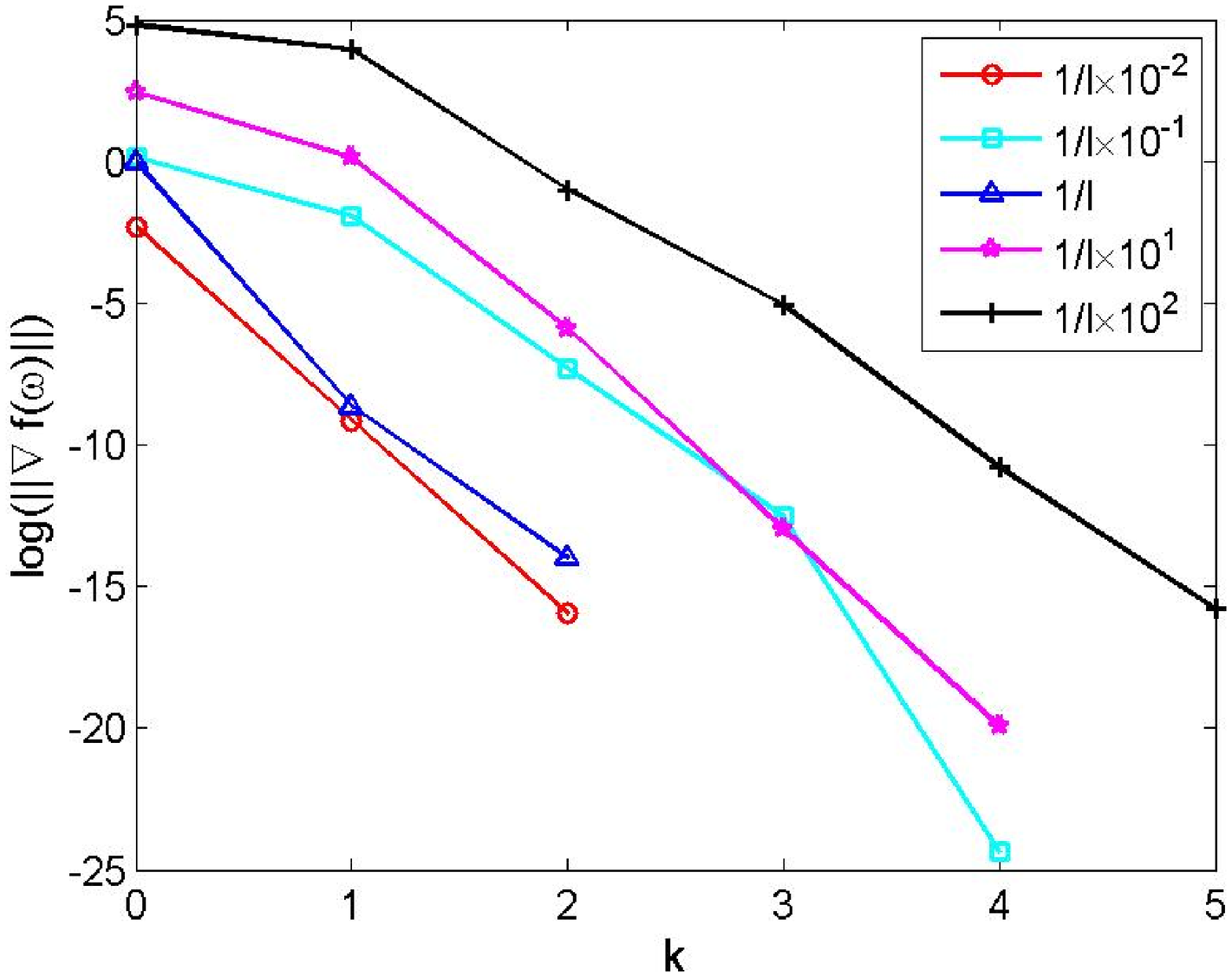}
  \caption*{(d)\ \ real-sim}
  \end{minipage}
  \caption{Performance Profile when $C\in \frac1l \times\{ 10^{-2}, 10^{-1}, 1, 10^{1}, 10^{2}\}$}\label{a}
\end{figure}

For L2-loss SVC, to show the quadratic convergence rate, we choose two data sets: ``w3a" and ``real-sim" to run Alg. \ref{algorithm2} and  plot the $\log\|\nabla f(\omega)\|$ during iterations  when $C\in \frac1l \times\{ 10^{-2}, 10^{-1}, 1, 10^{1}, 10^{2}\}$ in Figure $\ref{a}$. One can see that $\log\|\nabla f(\omega)\|$ decreased fast and stopped successfully within small iterations (the numbers of iterations in the two datasets are samller than 10). We can see that $\log\|\nabla f(\omega)\|$ decreases almost linearly along k, indicating the superlinear convergence rate of the semismooth Newton method.  

\begin{figure}[H]
  \centering
  \begin{minipage}[t]{.48\linewidth}
\includegraphics[width=1\textwidth]{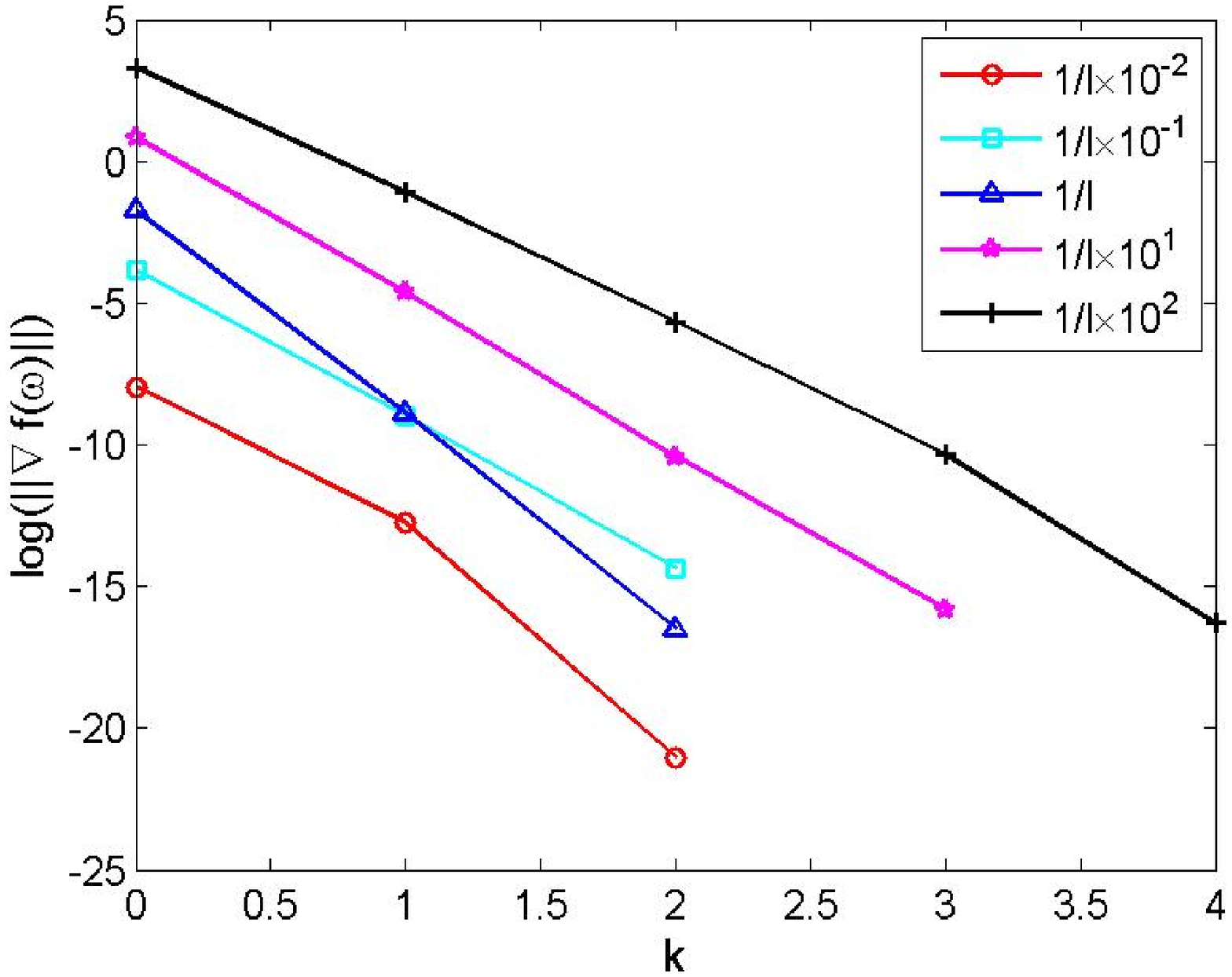}
  \caption*{(a) abalone}
  \end{minipage}
%
%
%
%
%
  \begin{minipage}[t]{.48\linewidth}

\includegraphics[width=1\textwidth]{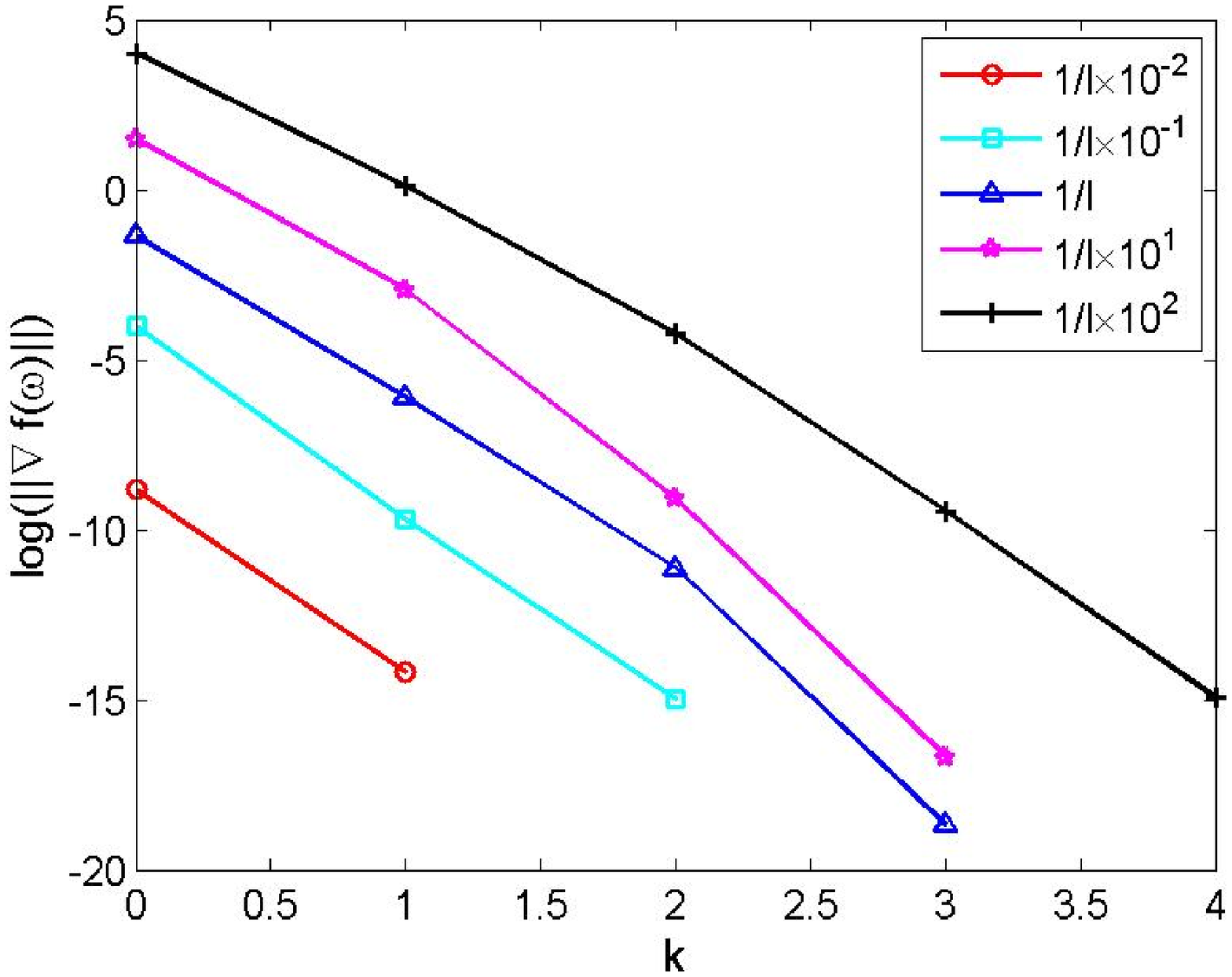}
  \caption*{(b) mpg}
  \end{minipage}
%
%
  \caption{$\log\|\nabla f(\omega)\|$ along Iterations with $C\in \frac1l \times\{ 10^{-2}, 10^{-1}, 1, 10^{1}, 10^{2}\}$ and $\epsilon =10^{-2}$}\label{b}
\end{figure}

 For $\epsilon$-L2 loss SVR, Figure $\ref{b}$ shows the trend of $\log\|\nabla f(\omega)\|$ during iterations via the semismooth Newton method in two data sets: ``abalone" and ``mpg".  We can observe that $\log\|\nabla f(\omega)\|$ decreases fast during iterations which again verifies the quadratic convergence rate of the semismooth Newton method.   
\subsection{Numerical Results for L2-loss SVC (\ref{prob})}
 In this part, to see the role of parameter $C$ in L2-loss SVC model (\ref{prob}), we test our algorithm with $C\in \frac1l \times\{ 10^{-2}, 10^{-1}, 1, 10^{1}, 10^{2}\}$ and report the results in  Table $\ref{B1}$ (We use each data set with 100$\%$ data). 
 \begin{table}[H]
 \centering
\caption{Numerical Results for L2-loss SVC with Different $C$'s.
}\label{B1}
\begin{tabular}{cccccc|cccccc}
\hline\noalign{\smallskip}
data&$C (\times \frac1l)$ & cg & k & res & t(s) &data&$C (\times \frac1l)$ & cg & k & res & t(s)\\
\noalign{\smallskip}\hline\noalign{\smallskip}
   \multirow {5}{*}{a1a} & $10^{-2}$ &5 & 3 & 1.80e-08 & 0.03   &\multirow {5}{*}{a2a}&  $10^{-2}$  & 5 & 3 & 1.80e-08 & 0.03 \\
   &$10^{-1}$ & 7 & 3 & 4.01e-07 & 0.03  && $10^{-1}$ & 7 & 3 & 4.02e-07 & 0.03 \\
    &$1$& 20 & 5 & 3.21e-08 & 0.07  && 1 &20 & 5 & 3.22e-08 & 0.07 \\
    &$10^{1}$& 47 & 7 & 7.69e-08 & 0.10   && $10^{1}$ &47& 7 &8.53e-08 & 0.10 \\
    & $10^{2}$ &121&9 & 1.87e-08 & 0.18&& $10^{2}$ & 123 & 9 &  2.05e-08& 0.18\\
   \noalign{\smallskip}\hline
   \multirow {5}{*}{a3a} & $10^{-2}$ &5 & 3 & 1.81e-08 & 0.03   &\multirow {5}{*}{a4a}&  $10^{-2}$  & 5 & 3 & 1.80e-08 & 0.03 \\
   &$10^{-1}$ & 7 & 3 & 4.01e-07 & 0.03  && $10^{-1}$ & 7 & 3 & 4.00e-07 & 0.03 \\
    &$1$& 20 & 5 & 3.31e-08 & 0.06  && 1 &20 & 5 & 3.32e-08 & 0.06 \\
    &$10^{1}$& 47 & 7 & 9.29e-08 & 0.09   && $10^{1}$ &47& 7 &6.46e-08 & 0.08 \\
    & $10^{2}$ &122 &9 & 3.49e-08 & 0.14&& $10^{2}$ & 121 & 9 & 2.98e-08& 0.14\\
   \noalign{\smallskip}\hline
   \multirow {5}{*}{a5a} & $10^{-2}$ &5 & 3 & 1.80e-08 & 0.03   &\multirow {5}{*}{a6a}&  $10^{-2}$  & 5 &3 & 1.80e-08 & 0.02 \\
   &$10^{-1}$ & 7 & 3 & 4.02e-07 & 0.03  && $10^{-1}$ & 7 & 3 & 3.97e-07 & 0.03 \\
    &$1$& 20 & 5 & 3.23e-08 & 0.05  && 1 &20 & 5 & 3.21e-08 & 0.04 \\
    &$10^{1}$& 47 & 7 & 6.60e-08 & 0.08  && $10^{1}$ &47& 7 &9.41e-08 & 0.06\\
    & $10^{2}$ &121 &9 & 2.88e-08 & 0.12&& $10^{2}$ & 123 & 9 & 1.93e-08& 0.10\\
   \noalign{\smallskip}\hline
   \end{tabular}
 \end{table}
 \begin{table}[H]
\begin{center}
 \begin{tabular}{cccccc|cccccc}
   \hline\noalign{\smallskip}
data&$C (\times \frac1l)$ & cg & k & res & t(s)  & data &$C (\times \frac1l)$ & cg & k & res & t(s)\\
\noalign{\smallskip}\hline\noalign{\smallskip}
   \multirow {5}{*}{a7a} & $10^{-2}$ &5 & 3 &  1.83e-08 & 0.02   &\multirow {5}{*}{a8a}&  $10^{-2}$  & 5 & 3 & 1.77e-08 & 0.02 \\
   &$10^{-1}$ &7 & 3 & 3.99e-07 & 0.02  && $10^{-1}$ & 7 & 3 & 3.97e-07 & 0.02 \\
    &$1$& 20 & 5 &  3.24e-08 & 0.04  && 1 &20 & 5 & 3.02e-08 & 0.04 \\
    &$10^{1}$& 47 & 7 & 9.89e-08 & 0.05   && $10^{1}$ &47& 7 &5.95e-08 & 0.06 \\
    & $10^{2}$ &119 &9 & 4.53e-08 & 0.07&& $10^{2}$ & 120 & 9 & 1.84e-08& 0.11\\
    \noalign{\smallskip}\hline
      \multirow {5}{1cm}{a9a} & $10^{-2}$ &5 & 3 & 1.80e-08 & 0.03   &\multirow {5}{1cm}{australian}&  $10^{-2}$  & 5 & 3 & 1.84e-08 & 0.00 \\
   &$10^{-1}$ & 7 & 3 & 4.02e-07 & 0.04  && $10^{-1}$ & 10 & 4 & 1.14e-08 & 0.00 \\
    &$1$& 20 & 5 &  3.23e-08 & 0.08  && 1 &23 & 5 & 3.85e-09 & 0.00 \\
    &$10^{1}$& 47 & 7 & 6.23e-08 & 0.11   && $10^{1}$ &45& 6 &3.57e-08 & 0.00\\
    & $10^{2}$ &121 &9 &2.47e-08 & 0.18&& $10^{2}$ & 73 & 7 & 2.02e-07& 0.01\\
   \noalign{\smallskip}\hline
    \multirow {5}{1cm}{breast-cancer}& $10^{-2}$ &5 & 4 & 1.32e-09 & 0.00   &\multirow {5}{1cm}{cod-rna}&  $10^{-2}$  & 14 & 6 & 2.28e-10 & 0.07 \\
   &$10^{-1}$ & 6 & 3 & 7.53e-07& 0.00  && $10^{-1}$ & 16 & 7 &4.82e-08 & 0.08\\
    &$1$& 13 &4 &  7.42e-07 & 0.00  && 1 &21 & 8 & 6.62e-09 & 0.10 \\
    &$10^{1}$& 34 & 8 & 4.10e-07 & 0.00   && $10^{1}$ &31& 11 &1.67e-08 & 0.13\\
    & $10^{2}$ &63 &9 & 1.17e-09 & 0.01&& $10^{2}$ & 44 & 13& 9.29e-09& 0.15\\
   \noalign{\smallskip}\hline
   \multirow {5}{1cm}{colon-cancer} & $10^{-2}$ &28 & 5 & 1.02e-08 & 0.02   &\multirow {5}{1cm}{diabetes}&  $10^{-2}$  & 4 & 3 &  2.10e-07 & 0.00 \\
   &$10^{-1}$ & 77 & 7 & 1.25e-08& 0.04  && $10^{-1}$ & 6 & 3 &3.09e-08 & 0.00 \\
    &$1$& 100 &7 & 1.88e-07 & 0.03  && 1 &11 & 4 & 2.84e-07 & 0.00 \\
    &$10^{1}$& 197& 12 & 7.44e-08 & 0.06   && $10^{1}$ &27& 6 &3.43e-09 & 0.00\\
    & $10^{2}$ &323 &19 & 7.50e-08 & 0.09&& $10^{2}$ & 30 & 6& 3.31e-07& 0.01\\
   \noalign{\smallskip}\hline
   \multirow {5}{1cm}{duke breast-cancer} & $10^{-2}$ &46 & 7 &6.46e-09& 0.06   &\multirow {5}{1cm}{fourclass}&  $10^{-2}$  & 3 & 3 &  1.66e-07 & 0.00 \\
   &$10^{-1}$ & 67 & 7 & 6.01e-07& 0.07  && $10^{-1}$ & 5 & 3 &9.09e-10 & 0.00 \\
    &$1$& 128 &11& 4.32e-08 & 0.11  && 1 &6 & 3 & 2.73e-08 & 0.00 \\
    &$10^{1}$& 207 & 18& 4.07e-08 & 0.18   && $10^{1}$ &12& 4 &7.03e-15 & 0.00\\
    & $10^{2}$ &406 &32 & 3.50e-07 & 0.36&& $10^{2}$ & 12 & 5& 7.38e-07& 0.00\\
   \noalign{\smallskip}\hline
   \multirow {5}{1cm}{german .numer} & $10^{-2}$ &6 & 3 &7.34e-10& 0.00   &\multirow {5}{1cm}{gisette}&  $10^{-2}$  & 18 & 5 &  2.54e-08 & 3.19 \\
   &$10^{-1}$ & 11 & 4 & 1.19e-08& 0.00  && $10^{-1}$ & 42 & 7 &3.74e-07 & 5.44 \\
    &$1$& 24 &5& 6.33e-09 & 0.00  && 1 &123 & 9 & 2.27e-08 & 9.00 \\
    &$10^{1}$& 48 &6 & 4.15e-07 & 0.01   && $10^{1}$ &292& 11 &9.16e-08 & 14.20\\
    & $10^{2}$ &87 &7 &1.42e-07 &0.01&& $10^{2}$ & 680 & 14& 6.55e-07& 24.13\\
   \noalign{\smallskip}\hline
   \multirow {5}{1cm}{heart} & $10^{-2}$ &5 & 3 & 1.18e-09& 0.00   &\multirow {5}{1cm}{ijcnn1}&  $10^{-2}$  & 4& 4 &  3.14e-08 & 0.06 \\
   &$10^{-1}$ & 10 & 4 & 4.59e-09& 0.00  && $10^{-1}$ & 6 & 3 &1.01e-08 & 0.05 \\
    &$1$& 18 &4& 7.14e-08 & 0.00  && 1 &11 & 4 & 7.53e-09 & 0.09 \\
    &$10^{1}$& 45 &6 & 3.35e-08 & 0.01   && $10^{1}$ &20& 5 &3.13e-07 & 0.13\\
    & $10^{2}$ &69 &7 &8.64e-08 &0.01&& $10^{2}$ & 49 & 7& 2.86e-07& 0.18\\
   \noalign{\smallskip}\hline
   \multirow {5}{1cm}{ionosphere} & $10^{-2}$ &5 & 3 & 6.00e-08& 0.00   &\multirow {5}{1cm}{leukemia}&  $10^{-2}$  & 39& 6 &  8.99e-07 & 0.05 \\
   &$10^{-1}$ & 10 & 4 &2.18e-08& 0.00  && $10^{-1}$ & 64 & 8 &9.30e-08 & 0.07 \\
    &$1$& 23 &5& 1.33e-08 & 0.00  && 1 &86 & 9 & 9.80e-07 & 0.08 \\
    &$10^{1}$& 46 &6 & 4.34e-07 & 0.01   && $10^{1}$ &282& 25 &9.25e-08 & 0.30\\
    & $10^{2}$ &99 &8 &4.61e-07 &0.01&& $10^{2}$ & 311 & 26& 9.25e-08& 0.28\\
   \noalign{\smallskip}\hline
    \multirow {5}{1cm}{liver-disorders} & $10^{-2}$ &23 & 5&7.05e-15& 0.00   &\multirow {5}{1cm}{mush rooms}&  $10^{-2}$  & 5& 3 & 1.47e-07 & 0.01 \\
   &$10^{-1}$ & 23& 5 &1.29e-07& 0.00  && $10^{-1}$ & 10 & 4 &4.84e-08 & 0.02 \\
    &$1$& 28 &7& 1.14e-11 & 0.00  && 1 &24 & 5 &6.03e-08 & 0.02 \\
    &$10^{1}$& 36 &8 &5.13e-13 & 0.00   && $10^{1}$ &56& 7 &5.52e-07& 0.03\\
    & $10^{2}$ &36 &8 &7.84e-12 &0.00&& $10^{2}$ & 145 & 10& 7.80e-08& 0.04\\
  \noalign{\smallskip}\hline
   \multirow {5}{1cm}{news20 .binary} & $10^{-2}$ &3 & 3&2.71e-07& 1.04  &\multirow {5}{1cm}{phishing}&  $10^{-2}$  & 3& 2 & 2.32e-07 & 0.01 \\
   &$10^{-1}$ & 5& 4 &4.46e-10& 1.44  && $10^{-1}$ & 5 & 3 &1.14e-08 & 0.02 \\
    &$1$& 7 &4& 6.22e-09 & 1.79  && 1 &8 & 4 &2.51e-08 & 0.03 \\
    &$10^{1}$& 12 &5  &3.68e-09 & 2.48   && $10^{1}$ &15& 5 &3.52e-07& 0.04\\
    & $10^{2}$ &19 &6 &2.04e-07 &3.60&& $10^{2}$ & 33 & 6& 9.53e-07& 0.05\\
  \noalign{\smallskip}\hline
   \multirow {5}{1cm}{rcv1 .binary} & $10^{-2}$ &3 & 3&1.59e-07& 0.10  &\multirow {5}{1cm}{real-sim}&  $10^{-2}$  & 3& 3 & 1.16e-07 & 0.20 \\
   &$10^{-1}$ & 6& 4 & 1.93e-10& 0.16  && $10^{-1}$ & 6 & 5 &2.59e-11& 0.32 \\
    &$1$& 7 &4&  1.17e-08 & 0.15  && 1 &5 & 3 &8.42e-07 & 0.24 \\
    &$10^{1}$& 12 &5  &6.07e-08 & 0.21   && $10^{1}$ &12& 5 &2.19e-09& 0.44\\
    & $10^{2}$ &24 &6 &2.31e-08 &0.33&& $10^{2}$ & 20 & 6& 1.38e-07& 0.60\\
   \noalign{\smallskip}\hline
   \end{tabular}\\
\end{center}
\end{table}
 \begin{table}[H]
\begin{center}
 \begin{tabular}{cccccc|cccccc}
   \hline\noalign{\smallskip}
data&$C (\times \frac1l)$ & cg & k & res & t(s)  & data &$C (\times \frac1l)$ & cg & k & res & t(s)\\
\noalign{\smallskip}\hline\noalign{\smallskip}
\multirow {5}{1cm}{skin $\_$nonskin} & $10^{-2}$ &15 & 5&4.18e-07& 0.18  &\multirow {5}{1cm}{splice}&  $10^{-2}$  & 8& 4 & 1.21e-09 & 0.01 \\
   &$10^{-1}$ & 18& 6 & 1.13e-08& 0.22  && $10^{-1}$ & 15 & 5 &3.11e-09& 0.01\\
    &$1$& 18 &6&  8.82e-08 & 0.22  && 1 &24 & 6 &3.84e-07 & 0.01 \\
    &$10^{1}$& 34 &11  &2.84e-07 & 0.43   && $10^{1}$ &38& 7 &1.08e-07& 0.02\\
    & $10^{2}$ &37 &12 &5.72e-08 &0.39&& $10^{2}$ & 54 & 8&1.22e-08& 0.02\\
   \noalign{\smallskip}\hline
   \multirow {5}{1cm}{sonar} & $10^{-2}$ &6& 3&4.68e-09& 0.00  &\multirow {5}{1cm}{svmguide1}&  $10^{-2}$  & 26& 8 & 5.70e-12 & 0.00 \\
   &$10^{-1}$ & 11& 4 & 1.71e-08& 0.00  && $10^{-1}$ & 29 &9 & 6.69e-09& 0.00\\
    &$1$& 23 &5 &  5.77e-07 & 0.00  && 1 &31 & 9 &7.89e-14 & 0.01 \\
    &$10^{1}$& 57 &7 &9.42e-08 & 0.01   && $10^{1}$ &32& 10 &1.09e-10& 0.01\\
    & $10^{2}$ &117 &8 &6.88e-07 &0.01&& $10^{2}$ & 39 & 11&4.10e-10& 0.01\\
  \noalign{\smallskip}\hline
\multirow {5}{1cm}{svmguide3} & $10^{-2}$ &4& 3&2.66e-07& 0.00  &\multirow {5}{1cm}{w1a}&  $10^{-2}$  & 4& 2 & 3.81e-07 & 0.04 \\
   &$10^{-1}$ & 6& 3 & 4.27e-08& 0.00  && $10^{-1}$ & 9 &4 & 2.82e-09& 0.07\\
    &$1$& 11 &4 &  7.32e-09 & 0.00  && 1 &14 & 4 &3.47e-07 & 0.08 \\
    &$10^{1}$& 24 &5 & 4.66e-07 & 0.01   && $10^{1}$ &43& 8 &5.04e-09& 0.14\\
    & $10^{2}$ &59 &7 & 7.58e-09 &0.01&& $10^{2}$ & 82 & 9&5.77e-07& 0.16\\
   \noalign{\smallskip}\hline
   \multirow {5}{1cm}{w2a} & $10^{-2}$ &4& 2&3.81e-07& 0.04  &\multirow {5}{1cm}{w3a}&  $10^{-2}$  & 4& 2 & 3.82e-07 & 0.04 \\
   &$10^{-1}$ & 9& 4 & 2.88e-09& 0.07  && $10^{-1}$ & 9 &4 & 2.89e-09& 0.07\\
    &$1$& 11 &4 &  3.31e-07 & 0.07  && 1 &14 & 4 &3.36e-07 & 0.07 \\
    &$10^{1}$& 43 &8 &4.60e-09 & 0.14   && $10^{1}$ &43& 8 &8.09e-09& 0.13\\
    & $10^{2}$ &81 &9 & 6.58e-07 &0.15&& $10^{2}$ & 81 & 9&6.03e-07& 0.15\\
   \noalign{\smallskip}\hline
   \multirow {5}{1cm}{w4a} & $10^{-2}$ &4 &2 & 3.76e-07& 0.03   &\multirow {5}{1cm}{w5a}&  $10^{-2}$  & 4& 2 &   3.85e-07 & 0.03 \\
   &$10^{-1}$ & 9 & 4 & 2.93e-09& 0.06  && $10^{-1}$ & 9 & 4 &2.99e-09 & 0.06 \\
    &$1$& 14 &4& 3.54e-07 & 0.06  && 1 &14 & 4 & 3.42e-07 & 0.06 \\
    &$10^{1}$& 43 &8 & 6.45e-09 & 0.12   && $10^{1}$ &43& 8 &1.42e-08 & 0.11\\
    & $10^{2}$ &81 &9 &5.83e-07 &0.14&& $10^{2}$ & 82& 9& 5.07e-07& 0.13\\
   \noalign{\smallskip}\hline
   \multirow {5}{1cm}{w6a} & $10^{-2}$ &4 &2 & 3.73e-07& 0.02   &\multirow {5}{1cm}{w7a}&  $10^{-2}$  & 4& 2 &   3.81e-07 & 0.02 \\
   &$10^{-1}$ & 9 & 4 & 2.83e-09& 0.05  && $10^{-1}$ & 8 & 4 &4.80e-08 & 0.04 \\
    &$1$& 14 &4& 3.49e-07 & 0.05  && 1 &14 & 4 & 3.50e-07 & 0.04 \\
    &$10^{1}$& 37 &7 & 9.03e-07 & 0.08   && $10^{1}$ &43& 8 &7.97e-09 & 0.07\\
    & $10^{2}$ &84 &9 &2.93e-07 &0.11&& $10^{2}$ & 84& 9&2.93e-07& 0.08\\
   \noalign{\smallskip}\hline
   \multirow {5}{1cm}{w8a} & $10^{-2}$ &4 &2 & 3.82e-07& 0.04   &\multirow {5}{1cm}{covtype .binary}&  $10^{-2}$  & 3& 3 &  1.58e-07 & 0.41 \\
   &$10^{-1}$ & 9 & 4 & 2.79e-09& 0.08  && $10^{-1}$ & 6 & 4 &7.49e-08 & 0.63 \\
    &$1$& 14 &4& 3.01e-07 & 0.07  && 1 &10 & 4 & 2.15e-07 & 0.73 \\
    &$10^{1}$& 43 &8 & 8.27e-09 & 0.15   && $10^{1}$ &20& 6 &7.64e-07 & 0.86\\
    & $10^{2}$ &82 &9 &6.44e-07 &0.16&& $10^{2}$ & 53& 11&4.84e-09& 1.16\\
   \noalign{\smallskip}\hline
   \end{tabular}\\
\end{center}
\end{table}
 From Table $\ref{B1}$, we summarize the following observations.
  \begin{enumerate}
    \item All the 210 tested instances are successfully solved by the semismooth Newton method. This suggests the semismooth Newton method is capable of solving problem (\ref{prob}) and the computation time of our algorithm is  small.
    \item  When our algorithm terminates, all the residuals (as shown in the column under ``res") are at least at the level of $10^{-6}$ within 10 iterations (Recall that the stopping criteria is $10^{-6}$), and even some residuals reached {$10^{-9}$ or $10^{-10}$. It indicates our algorithm  can stop successfully under the stopping criteria and return solutions of high accuracy}.  That is, the semismooth Newton method is effective to solve L2-loss SVC.
  \item In terms of different choices for  $C\in \frac1l \times\{ 10^{-2}, 10^{-1}, 1, 10^{1}, 10^{2}\}$, the semismooth Newton method can obtain the optimal solution even for different $C$'s. We can notice that the smaller the $C$, the less the number of iterations of our algorithm and the faster our algorithm can converge.
      \item  Our algorithm can converge to the optimal solution for most data sets within 10 iterations. 

  \end{enumerate}



  \subsection{Numerical Results  for $\epsilon$-L2-loss SVR (\ref{svr})}
   For $\epsilon$ in $\epsilon$-insensitive loss function \cite{Q15}, Ho et al. \cite{Q14}  performed experiments with and without using $\epsilon$ via the dual coordinate descent method. The  results indicate that for $\epsilon$-L2-loss SVR (\ref{svr}) , MSE is similar for different $\epsilon$'s. As a result, we fix $\epsilon=10^{-2}$ and test our algorithm with  different choices of $C$ for (\ref{svr}) since $\epsilon$ is insensitive. The results are reported in  Table $\ref{B}$.
 \begin{table}[H]
 \centering
\caption{Numerical Results for $\epsilon$-L2-loss SVR with Different $C$'s.}
\label{B}
\begin{tabular}{cccccc|cccccc}
\hline\noalign{\smallskip}
data&$C (\times \frac1l)$ & cg & k & res & t(s)  & data &$C (\times \frac1l)$ & cg & k & res & t(s)\\
\noalign{\smallskip}\hline\noalign{\smallskip}
    \multirow {5}{1cm}{abalone} & $10^{-2}$ &5 & 3 & 7.38e-10 & 0.00   &\multirow {5}{1cm}{bodyfat}&  $10^{-2}$  & 4 & 2 & 8.48e-08 & 0.00 \\
   &$10^{-1}$ & 6 & 3 & 5.74e-07 & 0.00  && $10^{-1}$ & 9 & 4 & 2.09e-10 & 0.00 \\
    &$1$& 11 & 3 & 6.95e-08 & 0.01  && 1 &14 & 4 & 6.67e-09 & 0.00 \\
    &$10^{1}$& 17 & 4 & 1.35e-07 & 0.01   && $10^{1}$ &27& 5 &8.83e-09 & 0.00 \\
    & $10^{2}$ &28&5 & 8.59e-08 & 0.01&& $10^{2}$ & 45 & 6 &  2.30e-07& 0.00\\
   \noalign{\smallskip}\hline
   \multirow {5}{*}{cpusmall} & $10^{-2}$ &4 & 2 & 1.79e-07 & 0.01   &\multirow {5}{1cm}{tfidf. train}&  $10^{-2}$  & 3 & 3 & 3.62e-07 & 1.71 \\
   &$10^{-1}$ & 7 & 3 & 5.76e-09 & 0.02  && $10^{-1}$ & 6 & 4 & 5.57e-09 & 2.87 \\
    &$1$& 12 & 4 & 2.77e-08 & 0.02  && 1 &6 & 3 & 2.07e-07 & 2.36 \\
    &$10^{1}$& 23 & 5 & 3.16e-08 & 0.03   && $10^{1}$ &10& 4 &4.06e-09 &3.58 \\
    & $10^{2}$ &43&5 & 2.19e-08 & 0.03&& $10^{2}$ & 13 & 4 &  2.08e-09& 3.87\\
   \noalign{\smallskip}\hline
   \multirow {5}{1cm}{tfidf. test} & $10^{-2}$ &3 & 3 & 1.16e-07 & 0.42   &\multirow {5}{1cm}{eunite 2001}&  $10^{-2}$  & 4 & 4 & 1.29e-07 & 0.00 \\
   &$10^{-1}$ & 6 & 4 & 3.72e-09 & 0.70  && $10^{-1}$ & 7 & 4 & 9.96e-08 & 0.00 \\
    &$1$& 6 & 3 & 9.14e-08 & 0.58  && 1 &10 & 4 & 7.28e-07 & 0.00 \\
    &$10^{1}$& 9 & 4 & 3.60e-08 & 0.92  && $10^{1}$ &21& 5 &8.52e-08 &0.00 \\
    & $10^{2}$ &12&4 & 4.97e-09 & 0.95&& $10^{2}$ & 49 & 6 & 9.08e-09& 0.00\\
   \noalign{\smallskip}\hline
   \multirow {5}{*}{housing} & $10^{-2}$ &6 & 3 & 1.47e-09 & 0.00   &\multirow {5}{*}{mg}&  $10^{-2}$  & 3 & 2 & 8.90e-07 & 0.00 \\
   &$10^{-1}$ & 8 & 3 & 6.49e-08 & 0.00  && $10^{-1}$ & 6 & 3 & 1.82e-08 & 0.00 \\
    &$1$& 17 & 4 & 8.42e-09 & 0.00  && 1 &10 & 3 & 2.83e-08 & 0.00 \\
    &$10^{1}$& 32 & 5 & 4.54e-08 & 0.01  && $10^{1}$ &15& 4 &1.93e-09 &0.00 \\
    & $10^{2}$ &40&5 & 3.30e-07 & 0.00&& $10^{2}$ & 18 & 4 & 1.49e-08& 0.00\\
   \noalign{\smallskip}\hline
   \multirow {5}{*}{mpg} & $10^{-2}$ &3 & 2 & 6.87e-07 & 0.00   &\multirow {5}{*}{pyrim}&  $10^{-2}$  & 6 & 3 & 2.98e-09 & 0.00 \\
   &$10^{-1}$ & 6 & 3 & 3.11e-07 & 0.00  && $10^{-1}$ & 10 & 4 & 5.27e-08 & 0.00 \\
    &$1$& 16 & 4 &7.66e-09 & 0.00  && 1 &22 & 5 & 4.56e-08 & 0.00 \\
    &$10^{1}$& 23 & 4 &5.68e-08 & 0.01  && $10^{1}$ &41& 5 &7.42e-07 &0.00 \\
    & $10^{2}$ &30&5 & 3.33e-07 & 0.01&& $10^{2}$ & 106 & 7 & 4.44e-07& 0.01\\
   \noalign{\smallskip}\hline
   \multirow {5}{1cm}{space$\_$ga} & $10^{-2}$ &3 & 3 & 6.09e-07 & 0.00   &\multirow {5}{*}{triazines}&  $10^{-2}$  & 8 & 4 & 1.77e-08& 0.00 \\
   &$10^{-1}$ & 6 & 3 & 2.13e-09 & 0.00  && $10^{-1}$ & 15 & 4 & 1.18e-08 & 0.00 \\
    &$1$& 9 & 3 &3.97e-07 & 0.00  && 1 &28 & 5 & 2.20e-08& 0.00 \\
    &$10^{1}$& 17 & 4 &2.47e-09 & 0.01  && $10^{1}$ &64& 6 &3.93e-08 &0.01 \\
    & $10^{2}$ &19&4 & 6.68e-09 & 0.01&& $10^{2}$ & 140 & 7 &2.81e-07& 0.01\\
   \noalign{\smallskip}\hline
    \end{tabular}
   \end{table}

     Table $\ref{B}$ shows that our algorithm can stop successfully under stopping criteria, which indicates the semismooth Newton method is  efficient for $\epsilon$-L2-loss SVR (\ref{svr}). All datasets are successfully solved by the semismooth Newton method in seconds. This suggests the semismooth Newton method is capable of solving problem (\ref{svr}) and the computation time of our algorithm is quite small. Our algorithm can converge to the optimal solution for all data sets within 7 iterations. Similarly, we can observe that the smaller the $C$, the less the number of iterations of our algorithm and the faster our algorithm can converge.

\subsection{Numerical Comparisions with LIBLINEAR}
 In this part, we compare our algorithm with some solvers in LIBLINEAR \footnote{We use the software LIBNIEAR version 2.11 downloaded from https://www.csie.ntu.edu.tw/~cjlin/liblinear/} which is the most popular and successful public software for support vector classification, regression and distribution estimation with linear kernel. We choose the following popular solvers for  L2-loss SVC and $\epsilon$-L2-loss SVR.
 \begin{itemize}
 \item DCD1 and TRON1: a dual coordinate descent method \cite{Q3} and a trust region Newton method \cite{Q4} for L2-loss SVC. 

\item TRON2 and DCD2: a trust region Newton method and a dual coordinate descent method \cite{Q14} for $\epsilon$-L2-loss SVR. 
\end{itemize}
We use a stratified selection to split each set to 60$\%$ training and 40$\%$ testing. For L2-Loss SVC, the training time and  accuracy (in percentage) on classification datasets are reported in Table \ref{C1} with {$C=\frac1l \times10^{2}$}. 

 \begin{table}[H]
 \centering
\caption{The Comparison Results  for L2-loss SVC. A1: DCD1; A2: TRON1; A3: {\color{red}the} Semismooth Newton Method}\label{C1}
\begin{tabular}{cccc}
\hline\noalign{\smallskip}
data&t(s) (A1$\mid$A2$\mid$A3) & accuracy (A1$\mid$A2$\mid$A3) \\
\noalign{\smallskip}\hline\noalign{\smallskip}
    a1a& 0.04$\mid$0.07$\mid$0.08&84.63$\mid$84.63$\mid$\textbf{84.66}\\
    a2a &0.03$\mid$0.06$\mid$0.08 &84.70$\mid$84.70$\mid$\textbf{84.72}\\
    a3a& 0.03$\mid$0.06$\mid$0.08&\textbf{84.67}$\mid$\textbf{84.67}$\mid$84.62\\
    a4a &0.03$\mid$0.05$\mid$0.08 &84.68$\mid$84.68$\mid$\textbf{84.73}\\
    a5a& 0.03$\mid$0.05$\mid$0.07&84.71$\mid$84.71$\mid$\textbf{84.74}\\
    a6a &0.02$\mid$0.03$\mid$0.05 &84.40$\mid$84.40$\mid$\textbf{84.95}\\
    a7a& 0.02$\mid$0.03$\mid$0.04&\textbf{84.78}$\mid$\textbf{84.78}$\mid$84.77\\
    a8a&0.03$\mid$0.04$\mid$0.06 &\textbf{84.31}$\mid$\textbf{84.31}$\mid$84.30\\
    a9a& 0.04$\mid$0.06$\mid$0.09&84.64$\mid$84.64$\mid$\textbf{84.66}\\
    australian &0.00$\mid$0.00$\mid$0.00 &84.78$\mid$84.78$\mid$\textbf{85.14}\\
    breast-cancer& 0.00$\mid$0.00$\mid$0.00& \textbf{98.90}$\mid$\textbf{98.90}$\mid$\textbf{98.90}\\
    cod-rna  &3.10$\mid$0.06$\mid$0.09 &81.58$\mid$\textbf{82.60}$\mid$76.01\\
    colon-cancer& 0.01$\mid$1.03$\mid$0.05& \textbf{72.00}$\mid$\textbf{72.00}$\mid$\textbf{72.00}\\     diabetes&0.00$\mid$0.00$\mid$0.00&\textbf{80.46}$\mid$\textbf{80.46}$\mid$79.48\\
    duke breast-cancer&0.02$\mid$1.75$\mid$0.17&\textbf{80.00}$\mid$\textbf{80.00}$\mid$\textbf{80.00}\\
    fourclass & 0.00$\mid$0.00$\mid$0.00&66.96$\mid$66.96$\mid$\textbf{74.94}\\
   german. numer & 0.01$\mid$0.00$\mid$0.01&76.50$\mid$76.50$\mid$\textbf{76.75}\\
     \noalign{\smallskip}\hline
     \end{tabular}
 \end{table}

 \begin{table}[H]
 \centering
 \begin{tabular}{cccc}
\hline\noalign{\smallskip}
data& t(s) (A1$\mid$A2$\mid$A3) & accuracy (A1$\mid$A2$\mid$A3) \\
\hline\noalign{\smallskip}
    gisette&4.93$\mid$12.12$\mid$14.18&\textbf{97.00}$\mid$\textbf{97.00}$\mid$\textbf{97.00}\\
    heart&0.01$\mid$0.00$\mid$0.01&85.19$\mid$85.19$\mid$\textbf{87.04}\\
    ijcnn1 & 0.08$\mid$0.07$\mid$0.08&91.44$\mid$91.44$\mid$\textbf{92.31}\\
    ionosphere& 0.01$\mid$0.00$\mid$0.01&\textbf{93.57}$\mid$\textbf{93.57}$\mid$92.86\\
    leukemia& 0.02$\mid$1.95$\mid$0.25& 26.67$\mid$26.67$\mid$\textbf{93.33}\\
    liver-disorders& 0.00$\mid$0.00$\mid$0.00&39.66$\mid$62.07$\mid$\textbf{65.52}\\
    mushrooms& 0.01$\mid$0.01$\mid$0.02& \textbf{96.43}$\mid$\textbf{96.43}$\mid$\textbf{96.43}\\
    news20.binary& 0.61$\mid$1.52$\mid$2.45&\textbf{72.14}$\mid$\textbf{72.14}$\mid$69.84\\
    phishing & 0.02$\mid$0.03$\mid$0.03& \textbf{90.59}$\mid$\textbf{90.59}$\mid$\textbf{90.59}\\
    rcv1.binary& 0.12$\mid$0.16$\mid$0.22&93.74$\mid$93.74$\mid$\textbf{94.07}\\
    real-sim&0.34$\mid$0.29$\mid$0.37& \textbf{78.78}$\mid$\textbf{78.78}$\mid$73.88\\
    skin$\_$nonskin & 15.78$\mid$0.08$\mid$0.17&89.16$\mid$89.16$\mid$\textbf{90.61}\\
    splice& 0.15$\mid$0.01$\mid$0.01& 84.94$\mid$84.94$\mid$\textbf{85.40}\\
    sonar& 0.00$\mid$0.00$\mid$0.01&14.46$\mid$14.46$\mid$\textbf{15.66}\\
    svmguide1 & 0.00$\mid$0.00$\mid$0.00& \textbf{11.89}$\mid$\textbf{11.89}$\mid$\textbf{11.89}\\
    svmguide3& 0.00$\mid$0.00$\mid$0.01&\textbf{40.44}$\mid$\textbf{40.44}$\mid$\textbf{40.44}\\
    w1a & 0.04$\mid$0.05$\mid$0.10& 99.32$\mid$99.32$\mid$\textbf{99.92}\\
    w2a& 0.05$\mid$0.06$\mid$0.09&99.31$\mid$99.31$\mid$\textbf{99.92}\\
    w3a & 0.05$\mid$0.06$\mid$0.09& 99.29$\mid$99.29$\mid$\textbf{99.93}\\
    w4a& 0.04$\mid$0.05$\mid$0.09&99.30$\mid$99.30$\mid$\textbf{99.92}\\
    w5a & 0.03$\mid$0.05$\mid$0.08& 99.27$\mid$99.27$\mid$\textbf{99.92}\\
    w6a& 0.03$\mid$0.03$\mid$0.06&99.36$\mid$99.36$\mid$\textbf{99.94}\\
    w7a & 0.02$\mid$0.02$\mid$0.05& 99.31$\mid$99.31$\mid$\textbf{99.95}\\
    w8a& 0.04$\mid$0.06$\mid$0.10&99.33$\mid$99.33$\mid$\textbf{99.91}\\
    covtype.binary& 31.25$\mid$1.18$\mid$0.70& 59.29$\mid$59.29$\mid$\textbf{61.54}\\
    \noalign{\smallskip}\hline
     \end{tabular}
\end{table}

From Table $\ref{C1}$, we get the following observations.
   \begin{enumerate}
     \item All of the three methods have high accuracy. The accuracy of most datasets was over 60$\%$, and even higher than 90$\%$ for some datasets. For the 42 classification data sets, compared with DCD1 and TRON1, the semismooth Newton method has same or higher classification accuracy for 34 datasets.

     \item The semismooth Newton method is competitive with DCD1 and TRON1 in  terms of cputime. In particular, for ``covtype.binary", the semismooth Newton method is much faster than DCD1 and TRON1. For ``skin$\_$nonskin", the semismooth Newton method takes shorter  time than DCD1 and is as fast as  TRON1.
     \end{enumerate}

 In  summary, the semismooth Newton method is very competitive with DCD1 and TRON1, in terms of accuracy and cputime.

 Next, we compare our algorithm with DCD2  and TRON2 for $\epsilon$-L2-loss SVR.
 The results are listed in Table \ref{C}.
  As we know, the smaller the MSE, the better the fitting of the model. {For $\epsilon$-L2-loss SVR,  we tested 12 regression data sets and we observed that when $C=\frac1l \times 10^{2}$ and $\epsilon=1e-2$, MSE via our algorithm is significantly smaller than DCD2 and TRON2 for all regression datasets.} 
   As for the cputime, these three methods are almost same. These indicate our algorithm is efficient and has better performance than DCD2 and TRON2.

 \begin{table}[H]
 \centering
\caption{The Comparison Results  for $\epsilon$-L2-loss SVR  with $C=\frac1l \times10^{2}$ and $\epsilon=1e-2$. B1: DCD2; B2: TRON2; B3: the Semismooth Newton Method}\label{C}
\begin{tabular}{cccc}
\hline\noalign{\smallskip}
 data& t(s) (B1$\mid$B2$\mid$B3) & MSE (B1$\mid$B2$\mid$B3)\\
\noalign{\smallskip}\hline\noalign{\smallskip}
    abalone& 0.00$\mid$0.00$\mid$0.01&50.07$\mid$50.07$\mid$\textbf{4.17}\\
    bodyfat & 0.00$\mid$0.00$\mid$0.00& 0.77$\mid$0.77$\mid$\textbf{0.00}\\
    cpusmall& 0.01$\mid$0.01$\mid$0.02&112.35$\mid$112.37$\mid$\textbf{102.24}\\
    tfidf.train & 1.64$\mid$1.43$\mid$1.73& 0.46$\mid$0.46$\mid$\textbf{0.14}\\
    tfidf.test& 0.57$\mid$0.41$\mid$0.78&0.40$\mid$0.40$\mid$\textbf{0.13}\\
    eunite2001 & 0.00$\mid$0.00$\mid$0.00& 131854 $\mid$131854 $\mid$\textbf{408.44}\\
    housing& 0.00$\mid$0.00$\mid$0.01&194.38$\mid$194.38$\mid$\textbf{71.45}\\
    mg& 0.00$\mid$0.00$\mid$0.00& 0.87$\mid$0.87$\mid$\textbf{0.02}\\
    mpg& 0.00$\mid$0.00$\mid$0.00&562.55$\mid$562.56$\mid$\textbf{37.48}\\
    pyrim & 0.01$\mid$0.00$\mid$0.01& 0.07$\mid$0.07$\mid$\textbf{0.01}\\
    space$\_$ga& 0.00$\mid$0.00$\mid$0.01&0.44$\mid$0.44$\mid$\textbf{0.03}\\
    triazines& 0.03$\mid$0.00$\mid$0.02&\textbf{0.03}$\mid$\textbf{0.03}$\mid$\textbf{0.03}\\
    \hline\noalign{\smallskip}
     \end{tabular}
 \end{table}

\subsection{Numerical Comparisons with SVRG-BB}
 In this part, we compare the semismooth Newton method with SVRG-BB  \cite{SVRG-BB} for the following squared hinge loss SVC:
 \be\label{svc}
\min_{\omega\in\Re^n} f_3(\omega):=\frac\lambda2\|\omega\|^2+\frac1l\sum_{i=1}^l \max(1-y_i(\omega^Tx_i), 0)^2.\\
\ee
We refer to \cite[Algorithm 3]{SVRG-BB} for the algorithm of SVRG-BB, and we use the following parameters in SVRG-BB: $m=2l,\ \omega^0=ones(n,1),\ \eta_0=0.1$.

Note that (\ref{svc}) and (\ref{prob}) are equivalent by choosing proper $C$ in (\ref{prob}). However,
when solving (\ref{prob}) by SVRG-BB, we find that SVRG-BB is sensitive to the selection of parameter $C$ in (\ref{prob}), and it  cannot converge for most data sets when $C\in \frac1l \times\{ 10^{-2}, 10^{-1}, 1, 10^{1}, 10^{2}\} $ or some other choices of $C$. Regarding this situation, we use the model $(\ref{svc})$ and we take $\lambda=1$, which is equivalent to our model $(\ref{prob})$ using $C=\frac1l$. Other settings are the same as before. Among the total 42 datasets, SVRG-BB cannot converge for some datasets (``cod-rna", ``colon-cancer", ``duke.breast-cancer", ``gisette", ``leukemia", ``liver-disorders", ``news20.binary", ``skin$\_$nonskin" and ``splice"). We test the rest 33 instances and we use a stratified selection to split each set to 60$\%$ training and 40$\%$ testing.

\begin{figure}[H]
  \centering
   \begin{minipage}[t]{.48\linewidth}
\includegraphics[width=1\textwidth]{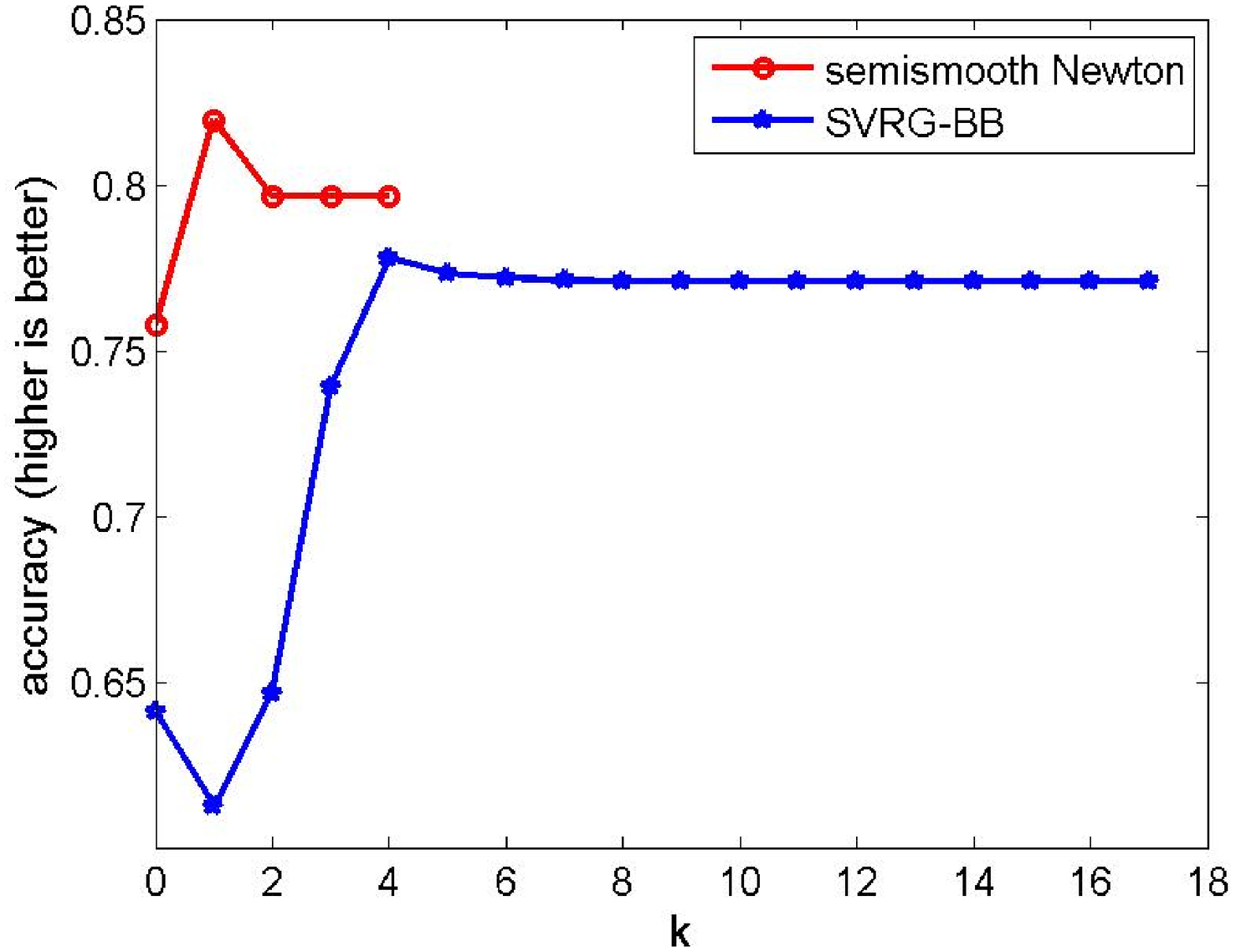}
  \caption*{(a)\ \  a5a}
  \end{minipage}
  \begin{minipage}[t]{.48\linewidth}
\includegraphics[width=1\textwidth]{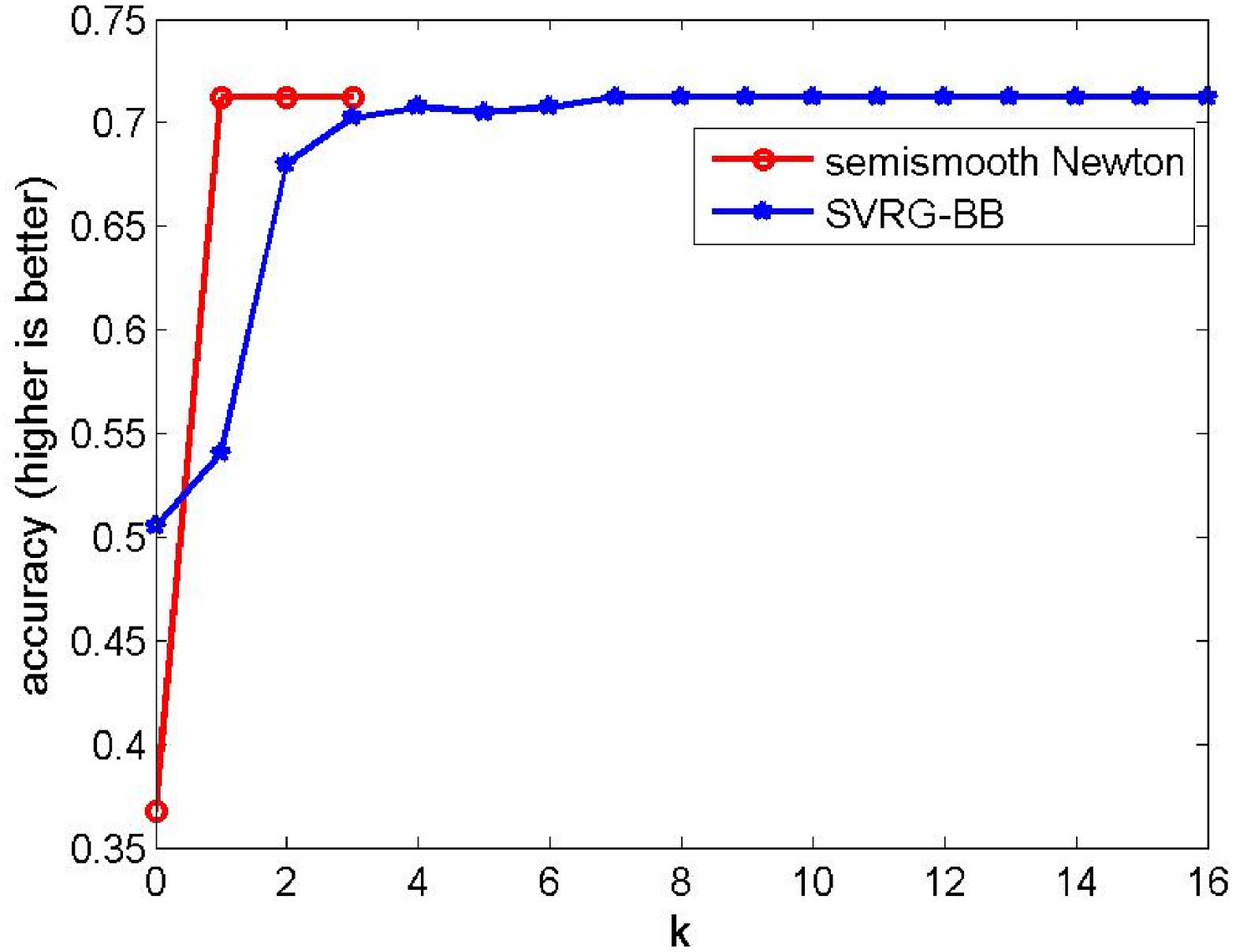}
  \caption*{(b)\ \  german.numer}
  \end{minipage}
  \begin{minipage}[t]{.48\linewidth}
\includegraphics[width=1\textwidth]{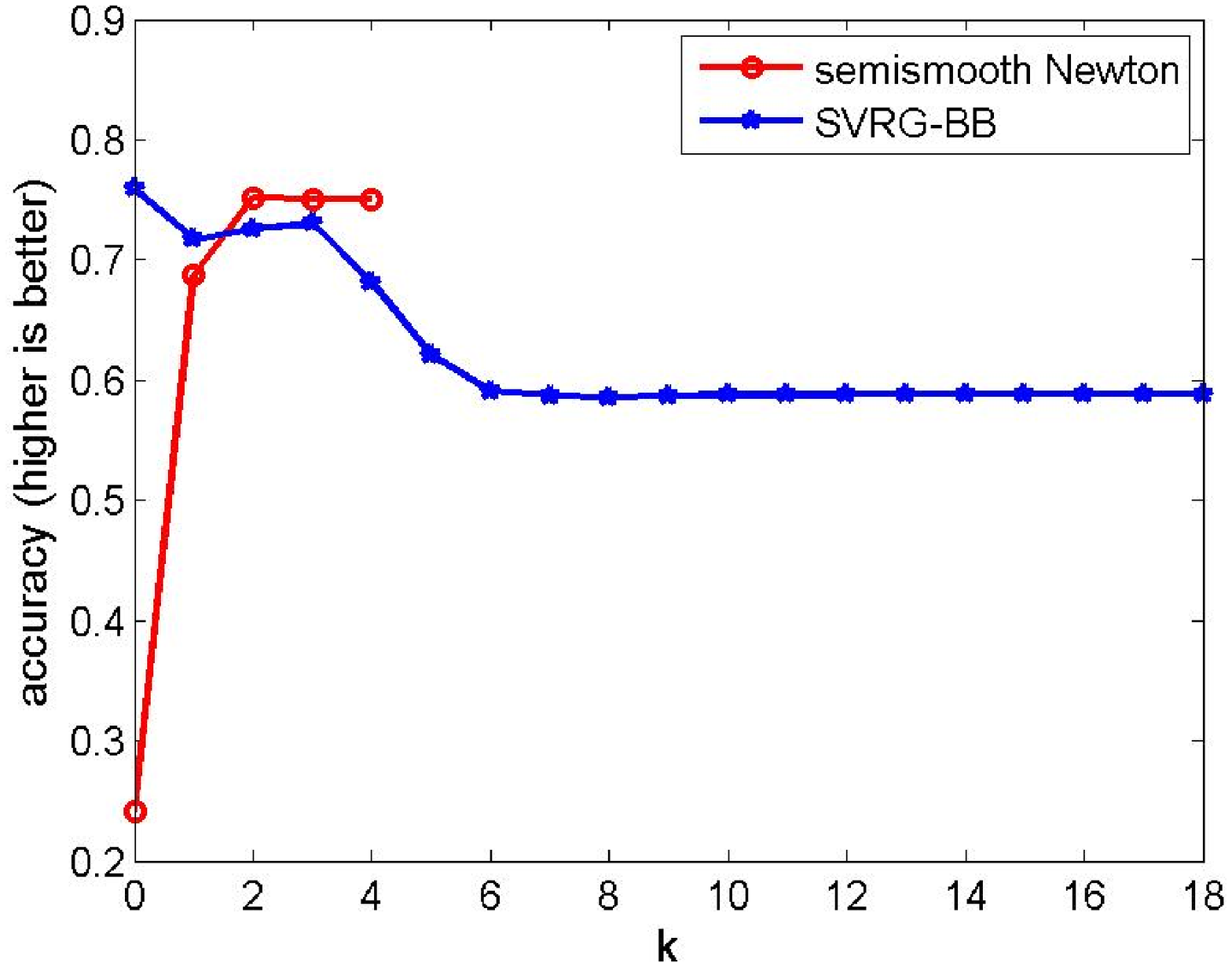}
  \caption*{(c)\ \ mushrooms}
  \end{minipage}
    \begin{minipage}[t]{.48\linewidth}
\includegraphics[width=1\textwidth]{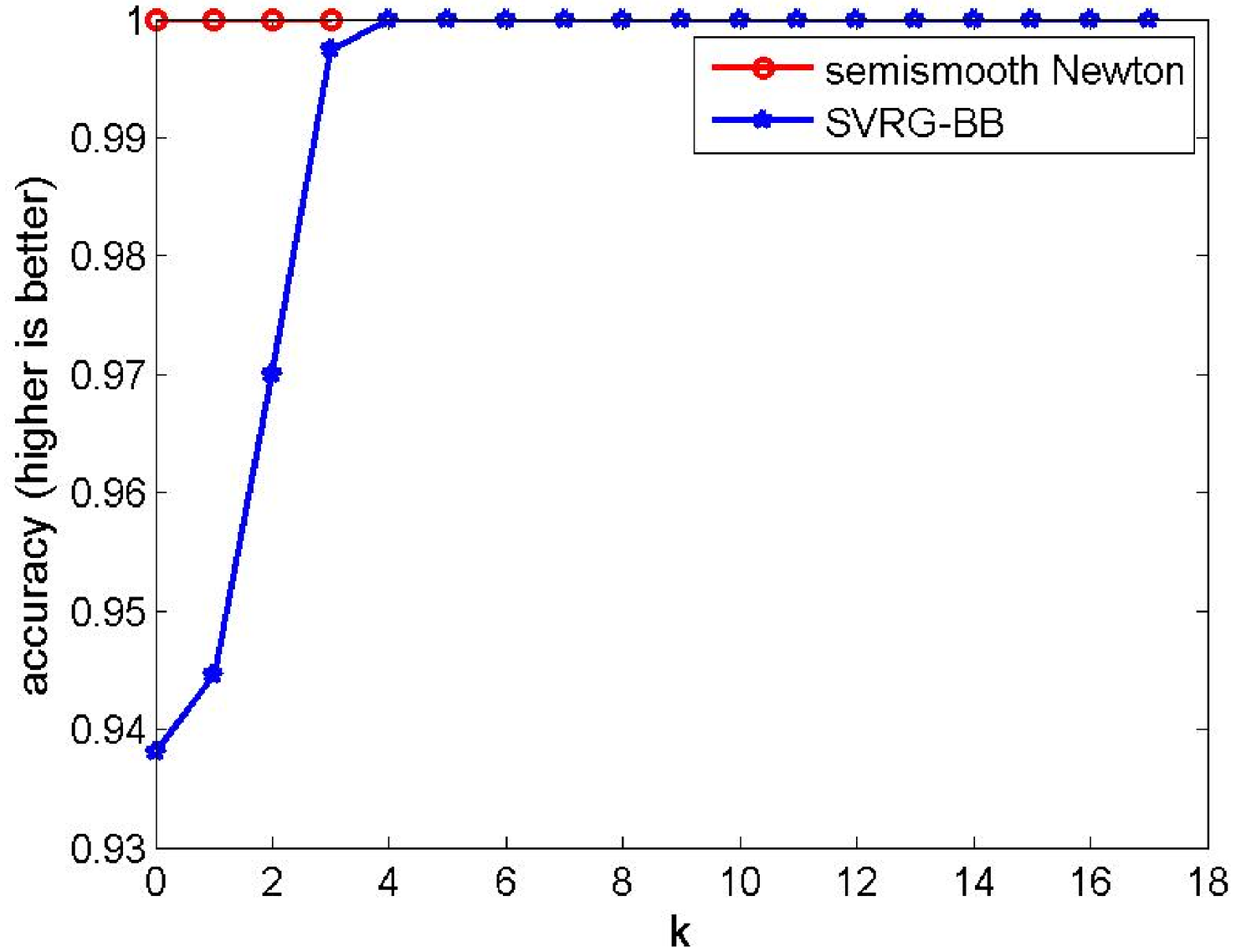}
  \caption*{(d)\ \ w5a}
  \end{minipage}
  \caption{Accuracy along Iterations of the Semismooth Newton Method and SVRG-BB  for (\ref{svc})}\label{R}
\end{figure}

In Figure \ref{R}, we selected 4 data sets: ``a5a", ``german.numer", ``mushrooms" and ``mushrooms" to show the accuracy along interations of the semismooth Newton Method and SVRG-BB  for (\ref{svc}). From Figure $\ref{R}$, we can see that our algorithm has much smaller iterations than SVRG-BB and the accuracy calculated by the semismooth Newton method is same or higher than SVRG-BB for these 4 datasets.

Next, we give the comparison results of the semismooth Newton method and SVRG-BB in Table \ref{D1}. From Table \ref{D1}, we have the following observations.
 \begin{enumerate}
     \item Our algorithm has smaller iterations than SVRG-BB. Our algorithm can satisfy the termination condition for most data sets within 5 iterations, however, SVRG-BB need to take about 20 iterations.
     \item It can be observed that the semismooth Newton method is significantly faster than SVRG-BB for all testing data sets. Our algorithm can converge to the optimal solution within a few seconds but SVRG-BB need to takes tens of seconds for most data sets. In particular, for datasets: ``rcv1.binary", ``real-sim" and ``covtype.binary", SVRG-BB takes 264, 484 and 632 seconds respectively, and our algorithm only  takes within 1 second.
     \item Both two methods have the high accuracy for most data sets. The accuracy via the semismooth Newton method is same or even higher than SVRG-BB for all data sets except ``ionosphere". For ``w1a" to ``w8a",  both two algorithms achieved 100$\%$ accuracy eventually.
     \end{enumerate}
 In summary, our algorithm is very effective and has better performance than SVRG-BB with regard to number of iterations, computational time and accuracy.

\begin{table}[H]
\centering
\caption{The Comparison Results of the Semismooth Newton Method and SVRG-BB.}\label{D1}
\begin{tabular}{ccccccc}
\hline\noalign{\smallskip}
  \multirow {2}{*}{data}&\multicolumn {3}{c}{semismooth Newton method} &\multicolumn {3}{c}{SVRG-BB}\\
  &k  & t(s) & accuracy &k  & t(s) & accuracy\\
\noalign{\smallskip}\hline\noalign{\smallskip}
a1a & 5 & 0.04 & \textbf{79.57} &18 & 38.53& 76.96\\
a2a & 5 & 0.04 & \textbf{79.57} &18 & 37.55& 77.03\\
a3a & 5 & 0.03 & \textbf{79.50} &19 & 36.54& 76.85\\
a4a & 5 & 0.03 & \textbf{79.60} &18 & 34.52& 77.00\\
a5a & 5 & 0.03 & \textbf{79.63} &18 & 30.88& 77.09\\
a6a & 5 & 0.02 & \textbf{79.80} &18 & 25.25& 77.16\\
a7a & 5 & 0.02 & \textbf{79.68} &18 & 19.62& 77.19\\
a8a & 5 & 0.03 & \textbf{79.35} &18 & 27.08& 76.45\\
a9a & 5 & 0.04 & \textbf{79.53} &18 & 38.63& 76.89\\
australian & 4 & 0.00 & \textbf{85.14} &16 & 0.71& 84.78\\
breast-cancer & 5 & 0.00 & \textbf{98.17} &15 & 1.98& \textbf{98.17}\\
diabetes & 4 & 0.00 & \textbf{72.96} &15 & 0.78& 69.71\\
fourclass & 4 & 0.00 & \textbf{74.78} &13 & 0.77& 71.01\\
german.numer & 4 & 0.01 & \textbf{71.25} &17 & 1.17& \textbf{71.25}\\
 heart & 4 & 0.00 & \textbf{84.26} &15 & 0.28& 83.33\\
 ijcnn1 & 4 & 0.09 & \textbf{90.37} &15 & 47.59& \textbf{90.37}\\
 ionosphere & 5 & 0.01 & 92.86 &17 & 0.43& \textbf{93.57}\\
 mushrooms & 5 & 0.02 & \textbf{75.08} &19 & 17.35& 58.77\\
  phishing & 4 & 0.02 & \textbf{57.73} &16 & 11.07& 55.88\\
  rcv1.binary & 4 & 0.11 & \textbf{51.93} &16 & 264.07& 51.91\\
  real-sim & 4 & 0.19 & \textbf{2.40} &16 & 483.65& 0.06\\
  sonar &6 & 0.00 & \textbf{7.23} &20 & 0.27& \textbf{7.23}\\
  svmguide1 &14 & 0.00 & \textbf{11.89} &34 & 4.48& \textbf{11.89}\\
  svmguide3 &4 & 0.00 & \textbf{40.44} &16& 1.30& \textbf{40.44}\\
  w1a &4 & 0.04 & \textbf{100} &18& 63.19& \textbf{100}\\
  w2a &4 & 0.04 & \textbf{100} &19& 62.30& \textbf{100}\\
  w3a &4 & 0.05 & \textbf{100} &19& 59.71&\textbf{100}\\
  w4a &4 & 0.04 & \textbf{100} &19& 57.71&\textbf{100}\\
  w5a &4 & 0.04 & \textbf{100} &18& 56.10& \textbf{100}\\
   w6a &4 & 0.03 &\textbf{100} &18& 41.32& \textbf{100}\\
   w7a &4 & 0.02 &\textbf{100} &19& 33.57& \textbf{100}\\
   w8a &4 & 0.05 &\textbf{100} &19& 66.85& \textbf{100}\\
   covtype.binary &4 & 0.43 & \textbf{61.54}&18& 631.54& \textbf{61.54}\\
   \hline\noalign{\smallskip}
\end{tabular}
\end{table}

\section{Conclusions}\label{sec-conclusions}
 In this paper, we apply the semismooth Newton method to solve two typical SVM models: the L2-loss SVC model and the $\epsilon$-L2-loss SVR model.  Our contribution in this paper is that by exploring the sparse structure of the models, we significantly bring down the computational complexity, meanwhile keeping the quadratic convergence rate. Extensive numerical experiments demonstrate the outstanding performance of the semismooth Newton method, especially for problems with huge size of sample data (for \texttt{news20.binary} problem with 19996 features and 1355191 samples, it only takes three seconds). In  particular, for the $\epsilon$-L2-loss SVR model, the semismooth Newton method significantly outperforms the leading solvers including DCD and TRON.
\end{document}